\documentclass[11pt]{article}

\pdfoutput=1

\usepackage{graphicx}
\usepackage{amsmath}
\usepackage{amssymb}
\usepackage{multirow}

\def\ee{\begin{equation}}
\def\eee{\end{equation}}
\def\bqn{\begin{eqnarray*}}
\def\eqn{\end{eqnarray*}}
\def\bnl{\begin{eqnarray}}
\def\enl{\end{eqnarray}}
\def\bma{\begin{bmatrix}}
\def\ema{\end{bmatrix}}
\def\bmx{\begin{matrix}}
\def\emx{\end{matrix}}
\def\ben{\begin{enumerate}}
\def\een{\end{enumerate}}
\def\bit{\begin{itemize}}
\def\eit{\end{itemize}}
\def\bei{\begin{itemize}}
\def\eei{\end{itemize}}
\def\bet{\begin{tabular}}
\def\eet{\end{tabular}}
\newcommand{\vn}[1]{\left|\left|#1\right|\right|}

\def\diag{{\rm diag}}

\def \col {\rm{col}}

\def\kN{{k=1,2,\ldots,N}}

\def\jn{{j=1,2,\ldots,n}}
\def \RR {\mathbb{R}}

\def\qedp{\hfill{$\blacksquare$}}
\def\qed{\hfill {$\square$}}
\newcommand{\refe}[1]{(\ref{#1})}
\def\salt{\vskip 0.4 true cm}
\newtheorem{theorem}{Theorem}
\newtheorem{corollary}{Corollary}
\newtheorem{definition}{Definition}

\newtheorem{remark}{Remark}

\newtheorem{lemma}{Lemma}

\date{}
\title{\bf Synchronization of Interconnected Systems\\ 
with Applications to Biochemical Networks:\\ 
an Input-Output Approach}
\author{L. Scardovi\thanks{Luca Scardovi is with the Department of Mechanical and Aerospace Engineering, Princeton University, USA. {scardovi@princeton.edu}. The work is supported in part by ONR grants N00014--02--1--0826 and N00014--04--1--0534.}\;, M. Arcak\thanks{Murat Arcak is with the Department of Electrical Engineering and Computer Sciences, University of California, Berkeley, USA. {arcak@eecs.berkeley.edu}.}\;, 
and E. D. Sontag\thanks{Eduardo Sontag is with the Department of Mathematics, Rutgers University, USA. {sontag@math.rutgers.edu}.  The work is supported in part by NSF grants 0504557 and 0614371, NIH grant 1R01GM086881, and AF grant FA9550.}
}
\begin{document}
\maketitle

\begin{abstract}
This paper provides synchronization conditions for networks of nonlinear
systems. The components of the network (referred to as ``compartments'' in this paper)
are made up of an identical interconnection of subsystems, each represented as
an operator in an extended $L_2$ space and referred to as a ``species''.  
The compartments are, in turn, coupled through a diffusion-like term among
the respective species.  The synchronization conditions are provided by combining the input-output properties of the
subsystems with information about the structure of network.  
The paper also explores results for state-space models, as well as biochemical applications. The work is motivated by cellular networks where signaling occurs both internally, through interactions of species, and externally, through
 intercellular signaling. The theory is illustrated providing synchronization conditions for networks of Goodwin oscillators.
\end{abstract}
\section{Introduction}
 The analysis of synchronization phenomena
 in networks
 has become an important topic in
 systems and control theory, motivated by
 diverse
 applications in physics,
 biology, and engineering. Emerging results in this area show that, in
 addition to the individual dynamics of the components, the network structure plays an important role in determining conditions leading to synchronization \cite{Hale,PhSl,StSe,Po,ScSaSe}.

In this paper, we study synchronization in networks of nonlinear systems, by
making use of the input-output properties of the subsystems comprising the
network.  Motivated by cellular networks where signaling occurs both
internally, through interactions of species, and externally, through
intercellular signaling, we assume that each component of the network
(referred to as a ``compartment'' in the paper) itself consists of subsystems
(referred to as ``species'') represented as operators in the extended $L_2$
space. The input to the operator includes the influence of other species
within the compartment as well as a diffusion-like coupling term between
identical species in different compartments.

A similar input-output approach was taken in
\cite{So1, Arcak:2006p424, Arcak:2008p782}
to study stability properties of individual compartments, rather than
synchronization of compartments.
These studies verify an appropriate \emph{passivity} property
\cite{sepulchre,vanderSchaft:2000p1141} for each species 
and form a ``dissipativity matrix'',
denoted here by $E$, that incorporates
information about the passivity of the subsystems, the interconnection
structure of the species, and the signs of the interconnection terms. 
To determine the stability of the network,
\cite{Arcak:2006p424, Arcak:2008p782} 
check the \emph{diagonal stability} of the dissipativity matrix, that is, the
existence of a diagonal solution $D>0$ to the Lyapunov equation $E^TD+DE<0$,
similarly to classical work on large-scale systems by Vidyasagar and others,
see \cite{vidyasagar2,sundareshan-vidyasagar,MoyHill78}.

In the special case of a {\it cyclic} interconnection structure with negative
feedback, this diagonal stability test encompasses the classical
\emph{secant criterion} \cite{TysOth78,Thr91}
used frequently in mathematical biology.
Following \cite{So1, Arcak:2006p424}, reference \cite{GuyBartStan:2007p3}
investigated synchronization of cyclic feedback structures using an
incremental variant of the passivity property.  This reference assumes that
only one of the species is subject to diffusion and modifies the secant
criterion to become a synchronization condition.

With respect to previous work, the main contributions of the present paper are
as follows: 
i) The results are obtained by using a purely input-output approach. This
approach requires in principle minimal knowledge of the physical laws
governing the systems, and is therefore particularly well-suited to
applications displaying high uncertainty on parameters and structure, such as
(molecular) biological systems. Results for systems with an ``internal
description'', i.e. in state space form, are derived as corollaries.
ii) The individual species are only required to satisfy an output-feedback passivity condition, compared to the stronger output-strict passivity condition in \cite{GuyBartStan:2007p3}.  
iii) The interconnections among the subsystems composing each network are not
limited to cyclic topologies, thus enlarging the class of systems for which
synchronization can be proved.
iv) The diffusive coupling can involve more than one species.
v) The new formulation allows exogeneous signals, and studies their effect on
synchronization.  

The paper is organized as follows.
In Section \ref{sec:ps} the notation used throughout the paper is summarized,
and the model under study is introduced. 
In Section \ref{sec:main} the main results are presented,
and the proofs can be found in Section \ref{sec:proof}. In Section \ref{sec:DiEx}
the operator property required to derive the synchronization condition is related to verifiable conditions for particular classes of systems described by ODE's; moreover the main results are extended to the case where the compartmental and the species couplings involve different variables. In Section \ref{sec:sec}, we show that the synchronization condition
can be expressed in terms of algebraic inequalities, for particular classes of
interconnection structures. Finally, in Section \ref{sec:go}, we illustrate the proposed theory, deriving synchronization conditions for a network of Goodwin oscillators.

\section{Preliminaries and problem statement}\label{sec:ps}
We denote by $L_{2e}$ the extended space of signals
$w : [0, \infty) \rightarrow \RR$ which have the property that each restriction $w_{T} = w|_{[0,T ]}$ is
in $L_{2}(0, T )$, for every $T > 0$. Given an element $w \in L_{2e}$ and any fixed $T > 0$, we
write $\vn{w}_T$  for the $L_2$ norm of the restriction $w_T$ , and given two functions
$v, w \in L_{2e}$ and any fixed $T > 0$, the inner product of $v_T$ and $w_T$ is denoted by
$\langle v, w\rangle_T$. The same notation is used for vector functions \footnote{We will denote by $L^{m}_{2e}$ the extended space of $m$ dimensional signals.}.

Consider $n$ identical \emph{compartments}, each composed of $N$ subsystems that we refer to as \emph{species}.
The input-output behavior of species $k$ in compartment $j$ is described by
\begin{equation}\label{ol}
y_{k,j} = H_{k} v_{k,j}, \quad k=1,\ldots,N, \quad j=1,\ldots,n,
\end{equation}
where $H_k$ is an operator to be further specified.
The interconnections among species and compartments is given  by:
\begin{equation}\label{control1}
v_{k,j} = w_{k,j} + \sum_{i=1}^{N} \sigma_{k,i}y_{i,j}  +\sum_{z=1}^{n}a^{k}_{j,z} (y_{k,z}-y_{k,j}), \quad k=1,\ldots,N, \quad j=1,\ldots,n,
\end{equation}
where the coefficients $ \sigma_{k,i} \in \RR,\,k,i =1,2,\ldots,N,$
represent the interconnection  between different species, and are identical in each
compartment. These coefficients are grouped into an $N\times N$  matrix:
\begin{equation}\label{sigmamat}\Sigma :=[\sigma_{k,i}], \quad k,i =1,2,\ldots,N,\end{equation}
and the resulting interconnection is called \emph{species
  coupling}.

The scalars $a^{k}_{j,z}, k=1,2\ldots,N,\, j,z=1,2\ldots,n$ are nonnegative
and represent the interconnection among systems of the same species in
different compartments. We will call this interconnection \emph{compartmental
  coupling}.
We assume that there are no self-loops, i.e. $a^{k}_{j,j}=0$,
$k=1,2\ldots,N,\, j=1,2\ldots,n$.
Note that different species
can possess different coupling structures (as implied by the superscript
$k$ in $a_{j,z}^k$).
The compartmental coupling is expressed in a diffusive-like form,
as a function of the differences between species in the respective
compartments, and not the species themselves. This is more general than
true diffusion, which would correspond to the special case in which
$a^{k}_{j,z} = a^{k}_{z,j}$ for all $k,j,z$;  under this symmetry condition,
the fluxes $a^{k}_{j,z} (y_{k,z}-y_{k,j})$ and $a^{k}_{z,j} (y_{k,j}-y_{k,z})$
(between the $k$th species in the $j$th and the $z$th compartments) would
cancel each other out.

\begin{figure}
 \centerline{
 \includegraphics[scale=0.42]{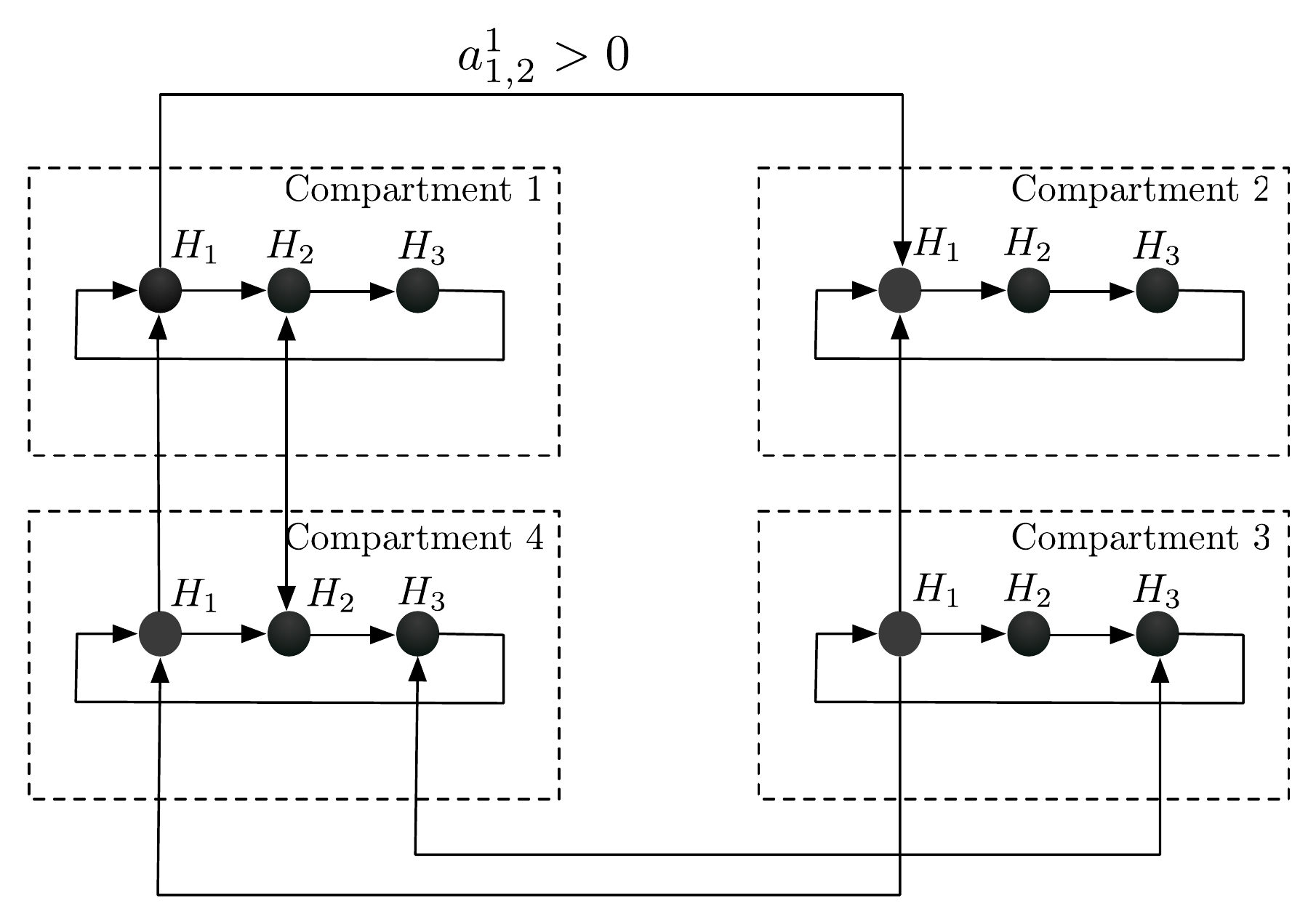}
 }
 \caption{Example of interconnection structure. Each compartment is composed by $3$ subsystems (represented as nodes of a graph) each characterized by an operator $H_k,\, k=1,2,3$. Two subsystems of the same species in different compartments are connected by an edge whenever the corresponding coefficient $a^k_{j,z}$ is positive. In each compartment, different species are interconnected according to a directed graph where the output of a system characterized by the operator $H_{i}$ enters as input of another system (characterized by an operator $H_{j}$) weighted by the coefficient $\sigma_{j,i}$.  In this example, the interconnections are cyclic ($\sigma_{i,j}=0$ unless $i=j+1$ mod $N$), but the theory allows arbitrary graphs. The compartments composing the network are assumed to be identical. For simplicity, no external inputs are shown in this figure.}
 \label{fig:compar2}
 \end{figure}

Finally, the scalars  $w_{k,j}$ are external inputs that can model e.g., $L_{2e}$ disturbances acting on the systems. The resulting
interconnected system can be represented as a graph as illustrated in
Figure \ref{fig:compar2}.

We denote by $Y_{k}=[y_{k,1},\ldots,y_{k,n}]^T$, $V_{k}=[v_{k,1},\ldots,v_{k,n}]^T$ and $W_{k}=[w_{k,1},\ldots,w_{k,n}]^T$  the vectors of the outputs, inputs and external signals of systems of the same species $k$, and  by $Y={\col}(Y_{1},\ldots,Y_N)$,  $V={\col}(V_{1},\ldots,V_N)$ and  $W={\col}(W_{1},\ldots,W_N)$  the stacked vectors.
We then rewrite the feedback law (\ref{control1}) as
\begin{equation}\label{control}
V_k (t) = W_{k}(t) + \sum_{i=1}^N \sigma_{ki} Y_i(t) -L_k Y_k(t) , \quad \quad k=1,2,\ldots,N,
\end{equation}
where $L_k, k=1,\ldots,N$ are Laplacian matrices associated to the compartmental coupling:
\[
l^k_{i,j}=\left\{
\begin{array}{ll}
\displaystyle \sum_{z=1}^{n} a^k_{i,z}, &   i=j\\
-a^k_{i,j},       &   i\ne j.
\end{array}
\right.
\]

The connectivity properties of the corresponding graphs are related to the algebraic properties of the Laplacian matrices and, in particular, to the notion of \emph{algebraic connectivity} extended to directed graphs in \cite{Wu:2005p964}:
\begin{definition}\label{connect}
For a directed graph with Laplacian matrix $L_k$, the \emph{algebraic
connectivity} is the real number defined as:
\begin{equation}\label{algebr}
\lambda_k = \min_{z\in {\cal P}} \frac{z^T L_k z}{z^Tz}
\end{equation}
where ${\cal P} = \{z\in \RR^n: z\, \bot\, 1_n, \vn{z} = 1\}$ and where $1_{n}\triangleq [1,1,\ldots,1]^{T} \in \RR^{n}$.
\end{definition}

To characterize synchronization mathematically, we denote the average of the outputs of the $n$ copies of the species
$k$ by:
\begin{equation}
\displaystyle \bar Y_{k} := \frac{1}{n}{1_n^T} Y_k,
\end{equation}
where
$1_{n}\triangleq [1,1,\ldots,1]^{T} \in \RR^{n}, \; \kN$, and define:
\begin{equation}\label{deviation}
\Delta Y_k := {\col}(y_{k1} - \bar Y_k, \ldots, y_{k,n} - \bar Y_k).
\end{equation}
Because $\Delta Y_{k}$ is equal to zero  if and
only if $Y_{k} = \alpha_{k} 1_{n}$  for some
$\alpha_{k}\geq 0$, $\vn{\Delta Y_{k}}_T$ measures the
synchrony of the outputs of the species  $k$ in the time
interval $[0,T]$.

We recall now an operator property that will be extensively used in the paper (the definitions are slightly adapted versions of those in \cite{CWillems:1971p1279}, \cite{Verma:2006p1287} and \cite{vanderSchaft:2000p1141}).
\begin{definition}\label{def:coco}
 Let $H: L^{m}_{2e} \rightarrow L^{m}_{2e}$. Then $H$ is relaxed cocoercive if there exists some $\gamma_c \in \RR$ such that for every pair of inputs $u,v \in  L^{m}_{2e}$
 \begin{equation}\label{coco}
\gamma_c \vn{Hu-Hv}_T^2 \leq \langle Hu-Hv , u - v\rangle_T , \quad \quad \forall T\geq0.
 \end{equation}
If \refe{coco} holds with $\gamma_c \geq 0$, then $H$ is called monotone. If \refe{coco} holds with $\gamma_c > 0$, then $H$ is called cocoercive.
\end{definition} 

\noindent Cocoercivity implies monotonicity and monotonicity
implies relaxed cocoercivity.  We refer to the maximum possible $\gamma_c$ with which \refe{coco} holds as the \emph{cocoercivity gain}, and
denote it as $\gamma$. The existence of $\gamma$ follows because the set of $\gamma_c$'s
that satisfy \refe{coco} is closed from above. In particular, we will call $\gamma$\emph{-relaxed cocoercive} the operators with a cocoercivity gain $\gamma \in \RR$. Notice that a $\gamma$-relaxed cocoercive operator with a strictly positive $\gamma$ is a cocoercive operator while, in general, it is only relaxed-cocoercive (monotone when $\gamma = 0$).

When there is no coupling between the compartments, i.e. $L_k = 0, \kN$, the
compartments are isolated and their stability depends on the
 species coupling. Stability with species coupling has been studied
in \cite{So1} with an input-output approach, and in
\cite{Arcak:2006p424,  Arcak:2008p782} with a Lyapunov approach.
Using the {\it output strict passivity} property
of the operators $H_k$:
 \begin{equation}\label{OSP}
 \gamma_k\|H_ku\|_T^2\le  \langle H_ku,u\rangle_T,
 \end{equation}
and defining the \emph{dissipativity matrix}\footnote{The matrix used in \cite{Arcak:2006p424, Arcak:2008p782} is slightly
different than the one used here; we are adopting the equivalent formulation
found in \cite{SoAr}.}
\begin{equation}\label{dissipmatrix}
E_{\gamma} =  \Sigma -\Gamma,\quad \quad \Gamma = \diag({\gamma_{1}},{\gamma_{2}},\ldots,{\gamma_{N}}), \quad \quad
\gamma=\col(\gamma_1,\cdots,\gamma_N)
\end{equation}
where $\Sigma$ is the interconnection matrix (\ref{sigmamat}),
\cite{Arcak:2006p424, Arcak:2008p782} prove stability of the interconnected system from the {\it diagonal stability} of the dissipative matrix $E_\gamma$; that is, from the existence of a diagonal matrix $D>0$ such that
\begin{equation}\label{diagstab}
E_{\gamma}^{T}D+DE_{\gamma} < 0.
\end{equation}

As we will see in Theorem 1 below,
the dissipativity matrix plays an
important role also when studying the synchronization properties of the system
\refe{ol}-\refe{control1}.  The key differences of Theorem 1 from \cite{Arcak:2006p424, Arcak:2008p782}  is that the output strict passivity property (\ref{OSP}) is replaced with the incremental property in Definition \ref{def:coco}, and the coefficients $\gamma_k$ in $E_\gamma$ are augmented with $\lambda_k$ terms from Definition \ref{connect}, which are due to diffusive coupling of the compartments.  Because the diagonal stability condition is more relaxed when $\gamma_k$ is augmented with $\lambda_k>0$, the new result makes it possible to show synchronization of the compartments when the individual compartments fail the stability test of \cite{Arcak:2006p424, Arcak:2008p782} and exhibit {\it e.g.} limit cycles.

\section{Main results}\label{sec:main}
The following theorem relates the properties of the interconnections and the operators to the synchrony of the outputs in the closed-loop system. In particular we show that, if the operators describing the open-loop systems are $\gamma$-relaxed cocoercive and the interconnection matrices satisfy certain algebraic conditions, the closed loop system has the property that external inputs with a ``high'' level of synchrony (as implied by a small $\|\Delta W\|_T$) produce outputs with the same property (small $\|\Delta W\|_T$).
\begin{theorem}
Consider the closed loop system defined by \refe{ol}-(\ref{control1}).
Suppose that the following assumptions are verified:
\begin{enumerate}
\item  Each operator $H_k$ is $\gamma_k$-relaxed cocoercive as in Definition \ref{def:coco}, $\kN$.
\item For $k=1,\ldots,N$, $\tilde \gamma_k := \lambda_{k}+\gamma_k>0$, where
  $\lambda_k$ is the algebraic connectivity in Definition \ref{connect} associated to
the matrix $L_k$ that describes the compartmental coupling of species $k$.
\item The dissipativity matrix $E_{\tilde \gamma}$ defined as in (\ref{dissipmatrix}) with
${\tilde \gamma} = {\col}(\tilde \gamma_1,\ldots,{\tilde \gamma_N})$, is
diagonally stable.
\end{enumerate}
Then, for all $w_{k,j}$, $y_{k,j}, \kN,\, \jn$ that satisfy \refe{ol} and \refe{control1} we have
\begin{equation}\label{ris}
\vn{\Delta Y}_{T} \leq \rho \vn{\Delta W}_{T},\quad \quad \forall T\geq0,
\end{equation}
for some $\rho>0$, and all $W \in L_{2e}^{Nn}$,  where
$\Delta W = {\col}(\Delta W_1,\ldots,\Delta W_N),
\Delta Y ={\col}(\Delta Y_1,\ldots,\Delta Y_N)$.
Moreover, if $W \in L_{2}^{Nn}$, then also $\Delta  Y \in L_{2}^{Nn}$, and
we have $\vn{\Delta  Y} \leq \rho \vn{\Delta W}$.
   \qed
\end{theorem}

Since subtracting a positive diagonal matrix from a diagonally stable matrix preserves diagonal stability, Theorem 1 says that the compartmental coupling increases the co-coercivity
 gain of a species whenever the algebraic connectivity is strictly positive. The algebraic connectivity is intimately related to topological properties of
 the underlying graph associated to the compartmental coupling (see Section 4,
 Remark 2).

This result can be extended to analyze synchronization in systems described with a state space formalism (with arbitrary initial conditions). This extension takes the form of a Corollary of Theorem 1. Consider the systems
\begin{equation}\label{statespace}
\begin{array}{rcl}
\dot x_{k,j} &=& f_{k}(x_{k,j},v_{k,j})\\
y_{k,j} &=& h_{k}(x_{k,j})
\end{array}
\quad k=1,\ldots,N, \quad j=1,\ldots,n,
\end{equation}
where $y_{k,j},u_{k,j}$ are scalars and $x_{k,j} \in \RR^p$ and the initial conditions are arbitrary.
We assume that $f_k(\cdot,\cdot)$ are locally Lipschitz in the first argument and that $h_k(\cdot)$ are continuous.
Furthermore  we assume that the systems are $L_2$-well-posed, in the sense that for each $v_{k,j} \in L_{2e}$ and each initial state there is a unique solution defined for all $t>0$ and the corresponding outputs $y_{k,j}(t)$ are also in $L_{2e}$.

If we set the initial conditions to zero we can define the input-output operators $H_{k}:L_{2e} \rightarrow L_{2e}$ by substituting any input $v_{k,j} \in L_{2e}$ in \refe{statespace}, solving the differential equation, and substituting the resulting state-space trajectory in $y_{k,j} = h_k(x_{k,j})$ in order to obtain the output function $y_{k,j}$. If we assume that the operators $H_k$ are well defined and we define the input as in \refe{control} then the results of Theorem 1 apply to the closed loop system.

The assumption that the initial state of the systems are set to zero is easy to dispose of, assuming appropriate reachability conditions. The following Corollary of Theorem 1 states that under the assumptions of Theorem 1 plus reachability and detectability conditions, the solutions of the compartments asymptotically synchronize. In particular we will assume that the closed-loop system \refe{statespace} is zero-state reachable, i.e. that for any state $x^*$ there exists an input belonging to $L_2$ that drives the system from the zero state to $x^*$ in finite time.

\begin{corollary}\label{coro}
Consider system (\ref{statespace}). Assume that the nonlinear operators associated to (\ref{statespace}) with zero initial conditions are well defined, that the conditions in Theorem 1 are verified and that the closed loop system is zero-state reachable. Then for all the outputs that satisfy \refe{statespace} with inputs as in (\ref{control1}) but with no external inputs ($w_{k,j}=0$), we have that  $\forall k =
1, \ldots, N, \forall i,j=1, \ldots, n,\; y_{k,i} (t)-y_{k,j}(t)\rightarrow 0$, as $t \rightarrow \infty$. In
addition, if for all initial states and all inputs any two state trajectories satisfy
\[
\vn{y_{k,j}-y_{k,i}}\rightarrow 0 \Rightarrow \vn{x_{k,j}-x_{k,i}} \rightarrow 0,\quad k=1,\ldots,N,\quad j,i\,=1,\ldots,n,
\]
as $t \rightarrow \infty$, then all bounded network solutions synchronize and the synchronized solution converges to the limit set of the isolated system (i.e. the system where $a_{k,j} = 0$ for every $k,j$).
\end{corollary} \qed

\section{Proof of the main result and corollary} \label{sec:proof}
We define the $(n-1)\times n$ matrix
\begin{equation}
Q=\left[ \begin{array}{ccccc}-1+(n-1)\nu & 1-\nu & -\nu & \cdots & -\nu \\
-1+(n-1)\nu & -\nu &1-\nu & \ddots & \vdots \\ \vdots & \vdots & \ddots & \ddots &-\nu \\
-1+(n-1)\nu & -\nu & \cdots & -\nu & 1-\nu
 \end{array}\right]
\end{equation}
where
\begin{equation}
\nu=\frac{n-\sqrt{n}}{n(n-1)}.
\end{equation}
It follows that $Q 1_n=0$, $QQ^T=I_{n-1}$, and
\begin{equation}\label{qtq}
Q^TQ=\left[\begin{array}{cccc} \frac{n-1}{n} & \frac{-1}{n} & \cdots & \frac{-1}{n} \\ \frac{-1}{n} & \frac{n-1}{n} & \ddots & \vdots \\ \vdots & \ddots & \ddots &\frac{-1}{n} \\ \frac{-1}{n} & \cdots & \frac{-1}{n} & \frac{n-1}{n} \end{array}\right] = I_n - \frac{1}{n}{1_n1_n^T}.
\end{equation}

By observing that
\begin{equation}
\tilde Y_{k}:= QY_{k}
 \end{equation}
 is equal to zero for every $k=1,\ldots,n$ if and only if $Y_{k} = \alpha_{k} 1_{n}$ for every $k=1,\ldots,N$, for some $\alpha_{k}\geq 0$, it is evident that also $\vn{\tilde Y_{k}}$ is a measure of synchrony for the outputs of the species (in different compartments) denoted by the index $k$. Moreover, since $Q^T Q Y_k =  \Delta Y_k$ from (\ref{deviation}) and (\ref{qtq}), $\tilde Y_k$ and $\Delta Y_k$ are related by $\Delta Y_k = Q^T \tilde Y_k$ and, thus,
 \begin{equation}
\vn{\Delta Y_k}_T^2 = \int_0^T \tilde Y_k^T Q Q^T \tilde Y_k\, dt = \vn{\tilde Y_k}_T^2, \quad \kN.
\end{equation}
In what follows we will use the same notation to measure input synchrony, i.e., we define $\tilde U_k = Q U_k$, $\tilde W_k = Q W_k$, $\tilde V_k = Q V_k$.
Before proving Theorem 1 we present a preliminary Lemma:
\begin{lemma}
Consider the open-loop systems (\ref{ol}). If the operators $H_k,\, k=1,2.\ldots,N$ are $\gamma_k$-relaxed cocoercive then
\begin{equation}\label{pass}
 \gamma_k \vn{\tilde Y_k}^2_T \leq \langle\tilde Y_{k},\tilde V_{k}\rangle_T,\quad \quad k=1,\ldots,N,
\end{equation}
for each $T>0$ and every $V_k \in L_{2e}^n$.
 \qed
\end{lemma}
\emph{Proof:} Consider the scalar product
\ee \label{1}
\langle \tilde V_{k}, \tilde Y_{k} \rangle_T
\eee
and define $z_{k,j} = v_{k,j} - \gamma_k y_{k,j}$ for every $k$, $j$, that in vector form reads
\ee \label{2}
Z_k = V_k - \gamma_k Y_{k}.
\eee
Define $\tilde Z_k = Q Z_k$. By substituting \refe{2} in \refe{1} we obtain
\ee \label{3}
\langle \tilde V_{k}, \tilde Y_{k} \rangle_T   = \langle \tilde Z_k , \tilde Y_{k} \rangle_T + \gamma_k \langle  \tilde Y_k, \tilde Y_{k} \rangle_T.
\eee
We first claim that the term $\langle \tilde Z_k , \tilde Y_{k} \rangle_T$ is nonnegative.
To show this, we use the $\gamma_k$-relaxed cocoercivity property of $H_k$ and obtain:
\ee \label{4}
\langle z_{k,i} - z_{k,j}, y_{k,i} - y_{k,j}\rangle_T = \langle v_{k,i} - v_{k,j}, y_{k,i} - y_{k,j}\rangle_T -  \gamma_k \langle y_{k,i} - y_{k,j}, y_{k,i} - y_{k,j}\rangle_T \geq 0,
\eee
for $i,j = 1,2,\ldots,n$.
By summing \refe{4} over $i,j=1,2,\ldots,n$ and by dividing by a normalization constant we get
\ee \label{5}
\frac{1}{2n}\sum_{i,j =1}^{n} \langle z_{k,i} - z_{k,j}, y_{k,i} - y_{k,j}\rangle_T = \langle Z_k, Y_k \rangle_T - n\langle \bar Z_k, \bar Y_k \rangle_T     \geq 0.
\eee
It follows that
\ee \label{6}
\langle \tilde Z_{k}, \tilde Y_{k} \rangle_T = \langle Z_{k},Q^T Q Y_{k} \rangle_T = \langle Z_{k},Y_{k} - {\bf 1}_n \bar Y_k \rangle_T =\langle Z_k, Y_k \rangle_T - n\langle \bar Z_k, \bar Y_k \rangle_T \geq 0,
\eee
which proves the claim.  Finally, from \refe{3} and \refe{6} we conclude that
\[
\langle \tilde V_{k}, \tilde Y_{k} \rangle_T   = \langle \tilde Z_k , \tilde Y_{k} \rangle_T + \gamma_k \langle  \tilde Y_k, \tilde Y_{k} \rangle_T \geq \gamma_k \vn{\tilde Y_k}_T^2,
\]
which is the desired inequality \refe{pass}.
\qedp \salt

\noindent We are now ready to prove Theorem 1.

\subsubsection*{Proof of Theorem 1}

Consider the inputs
\begin{equation}\label{input}
V_k (t) = U_k (t) -L_k Y_k(t) ,
\end{equation}
where $L_k$ are the Laplacian matrices representing the coupling between the compartments and the $U_k(t)$ are for now thought as external inputs.
From Lemma 1 and substituting (\ref{input}) in (\ref{pass}) we get,
\begin{equation}\label{first}
 \gamma_k \vn{\tilde Y_k}^2_T \leq  \langle\tilde Y_k,\tilde U_k\rangle_T -  \langle\tilde Y_k,Q L_{k} Y_k\rangle_T.
\end{equation}
Next, we note that $I_{n}-Q^TQ = \frac{1}{n} 1_n 1_n^T$ is a projection matrix onto the span of $1_n$. Because $L_k1_n=0$, it follows that $L_k(I_{n}-Q^TQ)Y_k=0$ and, thus,
\begin{equation}\label{nine}
L_k Y_k=L_kQ^TQY_k=L_kQ^T\tilde{Y}_k.
\end{equation}
By using (\ref{nine}) as well as the fact that
\[
Y_k^T(t)Q^TQL^TQ^TQY(t) = Y_k^T(t)Q^TQLQ^TQY(t)
\]
(because this expression is a scalar), we observe that:
\begin{equation}\label{second}
\langle\tilde Y_k,QL_{k} Y_{k}\rangle_T  = \frac{1}{2}\int_0^T \tilde Y_k^T(t)Q (L_k+L_k^T)Q^T \tilde Y_k(t) dt \geq \lambda_k\int_0^T \tilde Y_k^T(t) \tilde Y_k(t) dt = \lambda_k \vn{\tilde Y_{k}}_T^2,
\end{equation}
were $\lambda_{k}$ are the smallest eigenvalues of the symmetric part of the ``reduced Laplacian matrices'', i.e., of the matrices $(1/2)Q(L_{k} + L_{k}^{T})Q^{T}$. By using the properties of the matrix $Q$ it is straightforward to check that $\lambda_k$ is the algebraic connectivity as defined in Definition 1. Combining (\ref{first}) and (\ref{second}) we obtain
\[
\gamma_k \vn{\tilde Y_k}^2_T \leq \langle\tilde Y_k,\tilde U_k\rangle_T - \lambda_k \vn{\tilde Y_k}^2_T.
\]
From Assumption 2 we have that $\lambda_k>-\gamma_k$ for $k=1,2,\ldots,n$. We conclude that
\begin{equation}\label{final}
 \vn{\tilde Y_k}^2_T \leq  \frac{1}{\tilde \gamma_k}\langle\tilde Y_k,\tilde U_k\rangle_T, \quad \quad k=1,2,\ldots,N
\end{equation}
where $\tilde \gamma_k = {\gamma_{k}+\lambda_k}$.
The rest of the proof follows by the same argument as that used in the proof of ``Vidyasagar Lemma'' in~\cite{SoAr}, applied to the resulting input-output system $\tilde U_k \rightarrow \tilde Y_k$ and by using the condition (\ref{final}). Namely, we apply the feedback
\begin{equation}\label{control2}
U_k = W_k + \sum_{j=1}^N \sigma_{kj} Y_j, \quad \quad k=1,2,\ldots,N
\end{equation}
to the resulting system, where  $\sum_{j=1}^N \sigma_{kj} Y_j$ represents the interconnection  between the different species. By defining $U = {\col}(U_1,\ldots,U_N)$ (we apply this convention in general to denote vectors stacking), we rewrite (\ref{control2}) as
\begin{equation}\label{control3}
U = W + (\Sigma \otimes I_n) Y.
\end{equation}
We define
\[
E _{\tilde \gamma }\triangleq \Sigma - \Gamma_{\tilde \gamma},
\]
where
\[
\Gamma_{\tilde \gamma} \triangleq \diag (\tilde \gamma_1,\ldots,\tilde \gamma_N).
\]
From Assumption 3 the matrix $E_{\tilde \gamma}$ is diagonally stable i.e. there exist positive constants $d_i,\, i=1,\ldots,N$ such that
\begin{equation}\label{diagonal}
DE _{\tilde \gamma}+ E_{\tilde \gamma}^TD<0,
\end{equation}
and $D = \diag(d_1,\ldots,d_N)$.
Choose $\alpha>0$ such that $DE_{\tilde \gamma}+E_{\tilde \gamma}^TD<-2 \alpha I_N$ and observe that
\[
\langle Dz,E_{\tilde \gamma} z\rangle_T = \frac{1}{2}\int_0^T z^T(t) (DE_{\tilde \gamma}+E_{\tilde \gamma}^TD) z(t) dt \leq -\alpha \int_0^T z^T(t) z(t) dt = -\alpha \vn{z}_T^2.
\]
From (\ref{final}) we can write $\langle d_k \tilde Y_k, \tilde U_k - {\tilde \gamma_k} \tilde Y_k\rangle_T\geq 0 $ for $k=1,2,\ldots,N,$ and therefore
\[
\langle(D\otimes I_{n-1}) \tilde Y, \tilde U - (\Gamma_{\tilde \gamma} \otimes I_{n-1}) \tilde Y\rangle_T\geq 0,
\]
where $\tilde U = {\col}(\tilde U_1, \ldots, \tilde U_N)$ and $\tilde Y = {\col}(\tilde Y_1, \ldots, \tilde Y_N)$. Substituting $\tilde U =  \tilde W + (\Sigma \otimes I_{n-1}) \tilde Y$, where $\tilde W = {\col}(\tilde W_1, \ldots, \tilde W_N)$, we obtain
\[
\langle(D\otimes I_{n-1}) \tilde Y, \tilde W + (E_{\tilde \gamma} \otimes I_{n-1}) \tilde Y\rangle_T\geq 0,
\]
and using the Cauchy-Schwartz inequality we write
\[
\beta \vn{\tilde W}_T\vn{\tilde Y}_T\geq \langle(D\otimes I_n) \tilde Y, \tilde W\rangle_T \geq -\langle(D\otimes I_n) \tilde Y, (E_{\tilde \gamma} \otimes I_n) \tilde Y\rangle_T \geq \alpha \vn{\tilde Y}_T^{2},
\]
for some $\beta>0$. We conclude that
\[
\vn{\tilde Y}_T \leq \rho \vn{\tilde W}_T, \quad \forall T\geq 0,
\]
for any $W \in L_{2e}^{Nn}$, where $\rho = \beta/ \alpha$. As a direct consequence, if $W \in L_{2}^{Nn}$ then $\vn{\tilde Y} \leq \rho \vn{\tilde W}$. We conclude by observing that $\vn{ \tilde Y_k}_T  = \vn{ \Delta Y_k}_T$ and $\vn{ \tilde W_k}_T  = \vn{ \Delta W_k}_T$ for every $T\geq 0$. 

\qedp

To prove Corollary 1 we follow an argument similar to the one used in \cite{SoAr} to prove stability of interconnected systems.

\subsubsection*{Proof of Corollary 1}
Consider system (\ref{statespace}) where the initial conditions $x_{k,j}(0)=x_{k,j}^{0}$ are arbitrary, the inputs are
\begin{equation}\label{controlf}
v_{k,j} = \sum_{i=1}^{N} \sigma_{k,i}y_{i,j}  +\sum_{z=1}^{n}a^{k}_{j,z} (y_{k,z}-y_{k,j}),
\end{equation}
for $k=1,\ldots,N$, $j=1,\ldots,n$, and let $x_{k,j}$ be the solutions of the closed loop system.
Consider now system (\ref{statespace}) with initial conditions $z_{k,j}(0)=0$ and inputs $v_{k,j} + w_{k,j}$. From zero-reachability, there exist inputs  $\hat w_{k,j}: [0,T]\rightarrow \RR$ such that the solutions at time $T$ reaches the states $x_{k,j}^{0}$, i.e. $z_{k,j}(T) = x_{k,j}^{0}$ for every $k,j$. Consider now the input
\begin{equation}\label{constep}
w_{k,j}(t) =
\left \{
\begin{array}{rcl}
\hat w_{k,j}& t \in [0,T]\\
0 &t>T
\end{array}
\right.
\end{equation}
and let $z_{k,j}(\cdot)$ be the solution with initial state $z_{k,j}(0)=0$ and input $w_{k,j}(t)$ defined in (\ref{constep}). From causality we observe that $z_{k,j}(T) = x_{k,j}(0)$ and therefore $z_{k,j}(t+T) = x_{k,j}(t),\, t>T$, and therefore studying the steady state behavior of $z_{k,j}(\cdot)$ is equivalent to studying the steady state behavior of $x_{k,j}(\cdot)$. Consider the outputs $s_{k,j}(\cdot)$ associated to the solutions $z_{k,j}(\cdot)$ (with zero initial conditions and inputs $w_{k,j}$). From Theorem 1 we know that  $\vn{\Delta S} \leq \rho \vn{\Delta W}$, where $\Delta S = {\col}(\Delta S_1,\ldots,\Delta S_N)$, $\tilde W = {\col}(\Delta W_1,\ldots,\Delta W_N)$. Since each input $w_{k,j}$ is in $L_{2}$, $\Delta W$ is in $L^{nN}_{2}$ and we have that $\Delta S$ is in  $L^{nN}_{2}$ as well. Since the solutions $z_{k,j}$ are bounded, from continuity of $h_{k}(\cdot)$ we conclude that $S$ and therefore $\Delta S$ is absolutely continuous (see e.g., \cite{DeLoRySo}). From Barbalat's Lemma we conclude that $\Delta S \rightarrow 0$ for  $t\rightarrow \infty$ that implies that also $\Delta Y \rightarrow 0$ for  $t\rightarrow \infty$. This proves output synchronization. State synchronization directly follows from the additional property that for all initial states and all inputs, given two state trajectories we have that
\begin{equation}\label{id}
\vn{y_{k,j}-y_{k,i}}^{2}\rightarrow 0 \Rightarrow \vn{x_{k,j}-x_{k,i}}^{2}\rightarrow 0\quad 
\end{equation}
for $k=1,\ldots,N$,\; $j,i\,=1,\ldots,n$, as $t \rightarrow \infty$.
All the solutions that exist for all $t \geq 0$ converge to the set where for every $k =
1, \ldots, N$, and  $i,j=1, \ldots, n,\; x_{k,i} (t) =x_{k,j}(t)$, as $t \rightarrow \infty$. Since $L_k 1_n = 0$ for every $k=1,\ldots,N$, for all bounded network solutions, the synchronized
solution converges to the limit set of the isolated system (i.e. the system where the compartments are uncoupled). \qedp

\section{Discussion and extensions} \label{sec:DiEx}

\subsection{Conditions for relaxed cocoercivity}

The results presented require that the operators describing the input-output
relation of the isolated systems are relaxed cocoercive. This condition must in
general be checked case by case. Therefore, it is of interest to provide verifiable conditions, for particular
classes of systems described by ODE's, that imply relaxed cocoercivity of the
correspondent input-output operator.

Consider a one-dimensional system of the following form:
\begin{equation}\label{statespace3}
\begin{array}{rcl}
\dot x &=& -f (x) + u\\
y &=& x
\end{array}
\end{equation}
with zero initial conditions.
Suppose that $f$ is a function such that for every
$\sigma_1, \sigma_2 \in\RR$
\begin{equation}\label{lip2}
(\sigma_1-\sigma_2)(f(\sigma_1) -f(\sigma_2)) \geq \gamma \left(\sigma_1 - \sigma_2\right)^2,\quad  \gamma \in \RR.
\end{equation}
We now prove that the associated input-output operator from $u$ to $y$ is $\gamma$-relaxed cocoercive.
Fix the initial condition to zero. Consider two input functions $u_{1}$ and $u_{2}$ belonging to $L_{2e}$ and the corresponding outputs $x_{1}$ and $x_{2}$ (belonging to $L_{2e}$ as well). Then,
\[
\frac{1}{2} \frac{d}{dt}  (x_{2} -x_{1})^{2} = \left( f (x_{1})-f (x_{2}) + u_{2} - u_{1}\right) \left(x_{2}-x_{1} \right).
\]
Since \refe{lip2} holds, we conclude that
\[
\frac{1}{2} \frac{d}{dt}  (x_{2} -x_{1})^{2} \leq -{\gamma} \left( x_{2} - x_{1} \right)^2 + \left( u_{2} - u_{1}\right)\left(  x_{2}-x_{1}\right).
\]
Integrating both sides and assuming zero initial condition we finally obtain
\[
0 \leq
\frac{1}{2}(x_{2}(T) -x_{1}(T))^{2}
\leq
-{\gamma} \vn{x_{2} - x_{1}}^{2}_{T} + \langle u_{2} - u_{1},  x_{2}-x_{1}\rangle_{T},
\]
where the norm and the inner product are taken in the $L_{2e}$ spaces, showing that the associated input-output operator is ${\gamma}$-relaxed cocoercive. This result particularizes to linear time invariant system and e.g., for the system $\dot x = -a x +b u$, $y=x$, we obtain $\gamma = a/b$.

We conclude this section by characterizing a class of memoryless operators. First note that Definition \refe{def:coco} particularizes to the special case in which the operator $H$ is a nonlinear function $y=h(x)$ (static nonlinearity) and therefore it is possible to calculate the co-coercivity gain of a static nonlinearity by directly applying \refe{coco}.

We show now that a monotone increasing and Lipschitz continuous static nonlinearity $h(\cdot)$, with Lipschitz constant $(1/\xi)>0$, is a $\xi$-relaxed cocoercive operator (with positive $\xi$). From the Lipschitz condition we know that for every $\sigma_1, \sigma_2 \in \RR$
\begin{equation}\label{stat}
|h(\sigma_1) - h(\sigma_2)|\leq \frac{1}{\xi}|\sigma_1 -\sigma_2|.
\end{equation}
Multiplying both sides of \refe{stat} by $|h(\sigma_1) - h(\sigma_2)|$ we obtain
\[
(h(\sigma_1) - h(\sigma_2))^2\leq \frac{1}{\xi}|\sigma_1 -\sigma_2||h(\sigma_1) - h(\sigma_2)|.
\]
Since $\xi>0$ and $h$ is monotone increasing we conclude that
\[
\xi (h(\sigma_1) - h(\sigma_2))^2\leq (\sigma_1 -\sigma_2)(h(\sigma_1) - h(\sigma_2)),
\]
and we conclude that $h(\cdot)$ is $\xi$-relaxed cocoercive with $\xi>0$.

\subsection{State coupling versus output coupling}
The present work is motivated by synchronization in models of biochemical networks where the compartmental coupling represents the diffusion of reagent concentrations (states of the systems) through the compartments. The results presented in Section 1 assume that the species diffuse through the outputs that, in general, could be nonlinear functions of the species concentrations. In other words the compartmental and species couplings involve the same variables and this could be non realistic in the modeling of biological systems. In this section we generalize the results of the previous sections to handle this situation.

Consider the system
\begin{equation}\label{statespace2}
\begin{array}{rcl}
\dot x_{k,j} &=& f_{k}(x_{k,j}, v_{k,j}) \\
y_{k,j} &=& h_{k}(x_{k,j})
\end{array},
\quad k=1,\ldots,N, \quad j=1,\ldots,n,
\end{equation}
where
\begin{equation}\label{controldiff}
v_{k,j} = w_{k,j} + \sum_{i=1}^{N} \sigma_{k,i}y_{i,j}  +\sum_{z=1}^{n}a^{k}_{j,z} (x_{k,z}-x_{k,j}), \quad k=1,\ldots,N, \quad j=1,\ldots,n,
\end{equation}
that corresponds, in vector notation, to

\begin{equation}\label{controldiff2}
V_k (t) = W_{k}(t) + \sum_{j=1}^N \sigma_{kj} Y_j(t) -L_k X_k(t) , \quad \quad k=1,2,\ldots,N.
\end{equation}
Suppose that the solutions of \refe{statespace2} are defined for every input in $L_{2e}$, and that the respective output is in $L_{2e}$ as well. We fix the initial conditions to zero and we define the nonlinear operators $H_k$ associated to \refe{statespace2}. Suppose that each operator $H_k$ can be factorized as $H_k = T_{k}G_{k}$ where $T_k :v_{k,j} \rightarrow x_{k,j}$ is a nonlinear operator and $G_{k}:x_{k,j} \rightarrow y_{k,j}$ is a static nonlinearity associated to the functions $h_k(\cdot)$ from the real line to itself. 

Suppose that the operators $H_k$ and $G_k$ are $\gamma_k$-relaxed cocoercive and $\xi_k$-relaxed cocoercive respectively. Consider now the closed loop system. We follow the same lines of the proof of Theorem 1 exploiting the cocoercivity of the functions $h_k(\cdot)$ (associated to the operator $G_k$).
Consider the inputs
\begin{equation}\label{input_bis}
V_k (t) = U_k (t) -L_k X_k(t) ,
\end{equation}
where $U_k(t)$ are (for now) external inputs.
Since $H_k$ are $\gamma_k$-relaxed cocoercive, from Lemma 1, we get
\begin{equation}\label{first_bis}
\gamma_k  \vn{\tilde Y_k}^2_T \leq  \langle\tilde Y_k,\tilde U_k\rangle_T -  \langle\tilde Y_k,Q L_{k} X_k\rangle_T, \kN.
\end{equation}
Next we observe that
\begin{equation}\label{second_bis}
\langle\tilde Y_k,QL_{k} X_{k}\rangle_T  = \int_0^T \tilde Y_k^T(t)Q L_k Q^T \tilde X_k(t) dt.\end{equation}
Let's fix $k$ and define $s = h_k(\sigma)$ and $z(s) =\sigma -\xi_k s$. We observe that $z(s)$ is monotone increasing, in fact for every $s_1, s_2 \in \RR$
\[
(s_1-s_2)(z_1-z_2) =  (s_1-s_2)\left(\sigma_1-\sigma_2 - \xi_k(s_1 - s_2) \right)\geq 0,
\]
where the last inequality follows from the fact that $h_k(\cdot)$ is $\xi_k$-relaxed cocoercive.
By defining $Z_k(Y_k) = {\col}(z_k(y_{k,1}),\ldots,z_k(y_{k,n}))$, we use the identity $X_k = \xi_k  Y_k + Z_k$ to rewrite \refe{second_bis} as
\[
 \frac{\xi_k}{2} \int_0^T \tilde Y_k^T(t)Q (L_k + L_k^T)Q^T \tilde Y_k(t) dt +  \int_0^T Y_k^T Q^T Q  L_k Q^T Q Z_k(Y_k) dt.
\]
Suppose that $1_n^T L_k = 0$ then $Y_k^T Q^T Q  L_k Q^T Q Z_k(Y_k) = Y_k^T  L_k  Z_k(Y_k) $.
Since $L_k$ are doubly hyperdominant and $\tilde z_k(\cdot)$ are monotone increasing functions, from Theorem 3.7 in \cite{CWillems:1971p1279} we obtain
\[
 \int_0^T Y_k^T L_k Z_k dt \geq 0.
\]
We conclude that
\begin{equation}\label{third_bis}
\langle\tilde Y_k,QL_{k} X_{k}\rangle_T \geq \xi_k \lambda_k \int_0^T \tilde Y_k^T(t) \tilde Y_k(t) dt = \xi_k \lambda_k \vn{\tilde Y_{k}}_T^2.
\end{equation}
Combining (\ref{first_bis}) and (\ref{third_bis}) we obtain
\[
\gamma_k \vn{\tilde Y_k}^2_T \leq \langle\tilde Y_k,\tilde U_k\rangle_T - \xi_k \lambda_k \vn{\tilde Y_k}^2_T.
\]
Therefore, if  $\xi_k \lambda_k>-\gamma_k$ for $k=1,2,\ldots,n$, we conclude that
\begin{equation}\label{final_bis}
 \vn{\tilde Y_k}^2_T \leq  \frac{1}{\tilde \gamma_k}\langle\tilde Y_k,\tilde U_k\rangle_T, \quad \quad k=1,2,\ldots,N,
\end{equation}
where $\tilde \gamma_k = {\gamma_{k}+\xi_k \lambda_k}$.

The rest of the the derivation follows the same lines of the proof of Theorem 1 where $\tilde \gamma_k = \lambda_k + \gamma_k$ are redefined by $\tilde \gamma_k = {\gamma_{k}+\xi_k \lambda_k}$. This leads to the following result.

\begin{theorem}
Consider the closed loop system defined by \refe{statespace2} and \refe{controldiff2} with zero initial conditions.  Assume that the input-output operators $H_{j},\, j=1,2,\ldots,N$ associated to \refe{statespace2} are well-defined and suppose that the following assumptions are verified:
\begin{enumerate}
\item  The operators $H_k$ and the functions $h_k(\cdot)$ are, respectively, $\gamma_k$-relaxed cocoercive and $\xi_k$-relaxed cocoercive respectively for $\kN$.
\item The Laplacian matrices satisfy the condition $1_n^T L_k = 0$ and for $k=1,\ldots,N,$ $\tilde \gamma_k = {\gamma_{k}+\xi_k \lambda_k}$, where $\lambda_k$ is the algebraic connectivity associated to the compartmental coupling.
\item The matrix $E_{\tilde \gamma}$, where ${\tilde \gamma} = {\col}(\tilde \gamma_1,\ldots,{\tilde \gamma_N})$, is diagonally stable.
\end{enumerate}
Then,
\[
\vn{\Delta Y}_{T} \leq \mu \vn{\Delta W}_{T}, \quad \quad \forall T\geq 0,
\]
for some $\mu>0$ and $W \in L_{2e}^{Nn}$. Moreover, if $W \in L_{2}^{Nn}$, we have $\vn{\Delta Y} \leq \mu \vn{\Delta W}$.

 Furthermore, if the closed loop systems are zero-reachable, the closed loop system \refe{statespace2} and \refe{controldiff2} with no external inputs ($W=0$) and arbitrary initial conditions has the property that the outputs of the compartments synchronize, i.e. $\forall k =
1, \ldots, N, \forall i,j=1, \ldots, n,\; y_{k,i} (t) =y_{k,j}(t)$, as $t \rightarrow \infty$. In
addition, if for all initial states and all inputs any two state trajectories satisfy
\[
\vn{y_{k,j}-y_{k,i}}\rightarrow 0 \Rightarrow \vn{x_{k,j}-x_{k,i}} \rightarrow 0,\quad k=1,\ldots,N,\quad j,i\,=1,\ldots,n,
\]
as $t \rightarrow \infty$, then all bounded network solutions synchronize and the synchronized solution converges to the limit set of the isolated system (i.e. the system where $a_{k,j} = 0$ for every $k,j$). \qed
\end{theorem}

\begin{remark}
The algebraic connectivity and the properties of the Laplacian matrices can be related to properties of the underlying interconnection graph associated to the compartmental coupling.
The condition $1^T L_k = 0$ required by Theorem 2 is equivalent to assuming that the underlying graphs are balanced  (i.e. that for each vertex the sum of the weights of the edges entering in one vertex is equal to the sum of the weights of the edges exiting from the same vertex). Furthermore, since the resulting Laplacian matrices are doubly hyperdominant with zero excess the condition implies that $L_k +L_k^T$ are positive semidefinite and therefore the algebraic connectivity $\lambda_k \geq 0$. If furthermore we assume that the graph is connected, then the algebraic connectivity is guaranteed to be strictly positive (see e.g., \cite{Wu:2005p964}).
\end{remark}

\begin{remark}
The results presented in the paper can be used to analyze and design nonlinear observers. To see this, consider two identical compartments interconnected through a unidirectional compartmental coupling (for convention going from compartment two to compartment one).  The interconnection can involve one or more of the species. By interpreting the diffusive coupling terms as output injection, we can regard the first compartment as a nonlinear model and the second one as its state-observer. Thus Corollary \ref{coro} can be used to provide provable conditions for the observer error to converge to zero. The fact that our formulation allows for directed (non symmetric) graphs is here fundamental. In Section \ref{sec:go} we will illustrate this idea with a specific example.      
\end{remark}

\section{Special structures and synchronization conditions} \label{sec:sec}
Our results are based on the condition that
$E_{\tilde \gamma}$ be diagonally stable, which is related to both the compartmental coupling (through the algebraic connectivity) and the species coupling (through the interconnection matrix $\Sigma$).
In this Section, we analyze a number of interconnection structures and provide conditions for the matrix $E_{\tilde \gamma}$ to be diagonally stable. These conditions take the form of inequalities that link the algebraic connectivities of the compartmental coupling with the relaxed cocoercivity gains of the operators $H_k$.

\subsubsection*{Cyclic systems}
 \begin{figure}
 \centerline{
 \includegraphics[scale=0.65]{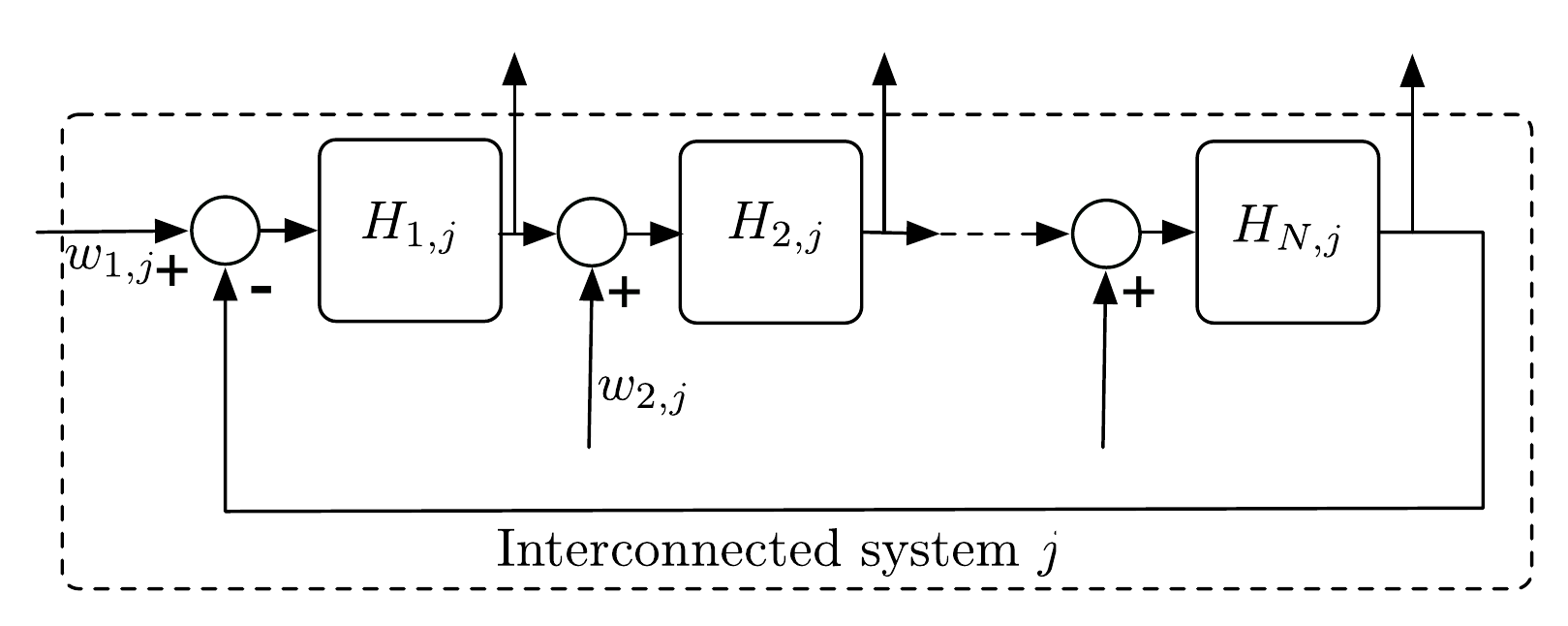}
 }
 \caption{Cyclic interconnection structure.}
 \label{fig:cyclic}
 \end{figure}
Stability of isolated cyclic systems is analyzed in \cite{Arcak:2006p424}.
Extending the work in \cite{Arcak:2006p424},
the output synchronization problem for cyclic systems was studied
in~\cite{GuyBartStan:2007p3}, for the special case in which the interconnected
systems are coupled \emph{through the output of the first system only}. Our approach is suitable for more general coupling (the interconnection
structure is depicted in Figure \ref{fig:cyclic}). The interconnection matrix is
 \[
 \Sigma_{\footnotesize \mbox{cyclic}}=
 \bma
 0 & 0 &\cdots& 0 & -1\\
 1 & 0 & 0&  \ddots & 0\\
 0 & 1 &0& \ddots&\vdots\\
 \vdots&\ddots&\ddots&\ddots& 0\\
 0 & \cdots &0& 1& 0
 \ema
 \]
 and the dissipativity matrix is therefore
 \[
 E_{\tilde \gamma}=
 \bma
- {\tilde \gamma_{1}}& 0 &\cdots& 0 & -1\\
 1 & -{\tilde \gamma_{2}}& 0&  \ddots & 0\\
 0 & 1 &-{\tilde \gamma_{3}}& \ddots&\vdots\\
 \vdots&\ddots&\ddots&\ddots& 0\\
 0 & \cdots &0& 1& -{\tilde \gamma_{N}}
 \ema.
 \]
 For this matrix to be diagonally stable the following secant condition must be satisfied \cite{Arcak:2006p424}:
 \[
 \frac{1}{\tilde \gamma_{1}} \frac{1}{\tilde \gamma_{2}}\cdots \frac{1}{\tilde \gamma_{N}}<\sec(\pi/N)^{N}.
 \]
 Since $\tilde \gamma_{k} =  {\gamma_{k}+\lambda_{k}}>0$, the secant condition leads to:
 \begin{equation}\label{sec_tot}
 \prod_{k=1}^N  \frac{1}{\gamma_{k}+\lambda_{k}} < \sec(\pi/N)^{N}.
 \end{equation}

Our approach generalizes the result of \cite{GuyBartStan:2007p3} (note that
$\gamma_k$ correspond to $1/\gamma_k$ in \cite{GuyBartStan:2007p3}) where the
coupling among the systems is limited to the first system (i.e., when
$\lambda_j =0$, $j=2,\ldots,N$ in (\ref{sec_tot})). In fact, in this case
\refe{sec_tot} reduces to
\begin{equation}\label{part}
\lambda_1 > \frac{\cos(\pi/N)^{N}}{\gamma_{2}\cdots \gamma_{N}}  - \gamma_1,
\end{equation}
which is the expression provided in \cite{GuyBartStan:2007p3}.

As an example, consider the case where each species in a compartment is
 directly connected to the respective species in each other compartment with the same weight $q$,
i.e. $a^k_{ij} = q$ for every $i,j =1,\ldots,n$ and $k = 1,\ldots N$. This
implies that the Laplacian matrices are
\[
L_k = q\, n \left(I_n -\frac{1}{n}1_n 1_n^T\right), \quad k=1,2,\ldots, N,
\]
and that $\lambda_k = q\,n$, $k=1,\ldots,N$. In this case, \refe{sec_tot}
specializes to:
\[
 \prod_{k=1}^N  \frac{1}{\gamma_{k}+q\, n } < \sec(\pi/N)^{N},
 \]
where, since $\lambda_k+\gamma_k$ must be strictly positive,  the condition
$q> -\gamma_k/n$, $k=1,2,\ldots,N$ must be satisfied. If we restrict the
compartmental coupling to only the first species, \refe{part} takes the simple
form
\[
q\,n >   \frac{\cos(\pi/N)^{N}}{\gamma_{2}\cdots \gamma_{N}} - \gamma_1.
\]

\subsubsection*{Branched structures}

  \begin{figure}
 \centerline{
 \includegraphics[scale=0.65]{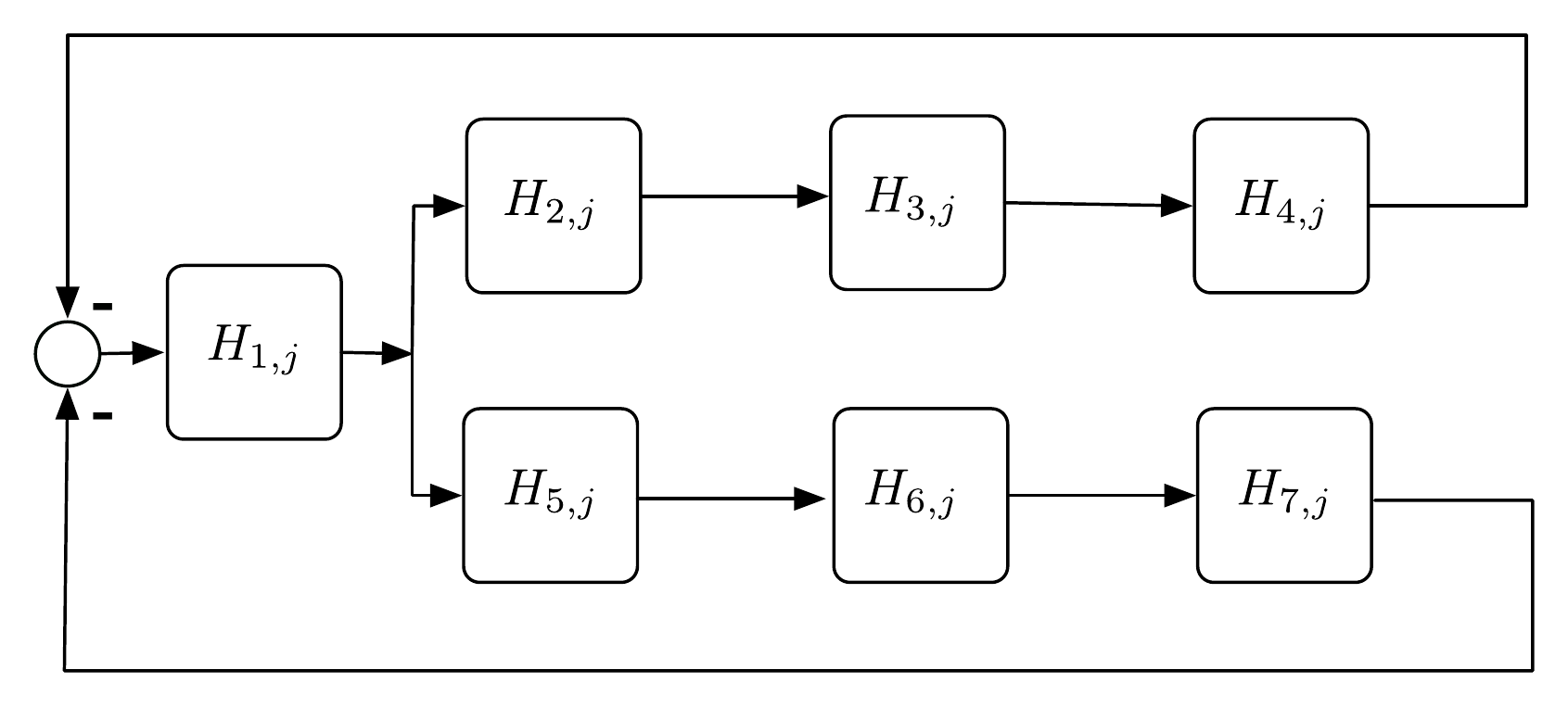}
 }
  \caption{Branched interconnection structure, for notational simplicity the inputs $w_{jk}$ are not shown in the picture.}
  \label{fig:branched}
 \end{figure}
In  \cite{Arcak:2008p782} several interconnection structures have been analyzed and diagonal stability is proven for the associated dissipative matrices. \\

\noindent i) For the interconnection structure depicted in Figure \ref{fig:branched}, the interconnection and dissipativity matrices are, respectively,
  \[
 \Sigma_{\footnotesize \mbox{b1}}=
 \bma
 0 & 0 & 0 & -1 & 0&0& -1\\
 1 & 0 & 0 & 0 & 0&0& 0\\
 0 & 1 & 0 & 0 & 0&0& 0\\
 0 & 0 & 1 & 0 & 0&0& 0\\
 1 & 0 & 0 & 0 & 0&0& 0\\
 0 & 0 & 0 & 0 & 1&0& 0\\
 0 & 0 & 0 & 0 & 0&1& 0\\
 \ema,
 \]
  \[
 E^{\footnotesize \mbox{b1}}_{\tilde \gamma}=
 \bma
  -{\tilde \gamma_1} & 0 & 0 & -1 & 0&0& -1\\
1 & -\tilde \gamma_{2}  & 0 & 0 & 0&0& 0\\
 0 & 1 & -{\tilde \gamma_{3} } & 0 & 0&0& 0\\
 0 & 0 & 1 & -{\tilde \gamma_{4} } & 0&0& 0\\
 1 & 0 & 0 & 0 & -{\tilde \gamma_{5} }&0& 0\\
 0 & 0 & 0 & 0 & 1&-{\tilde \gamma_{6} }& 0\\
 0 & 0 & 0 & 0 & 0&1& -{\tilde \gamma_{7} }\\
 \ema.
 \]
 Lemma 2 in \cite{SoAr} shows that $ E^{\footnotesize \mbox{b1}}_{\tilde \gamma}$ is diagonally stable iff the condition:
 \[
\frac{1}{\tilde \gamma_{1} \tilde\gamma_{2} \tilde \gamma_{3} \tilde \gamma_{4}} +\frac{1}{ \tilde\gamma_{1} \tilde\gamma_{5} \tilde\gamma_{6} \tilde\gamma_{7}} <\sec(\pi/4)^{4}  \]
  holds. Since $\tilde \gamma_{k} =  {\gamma_{k}+\lambda_{k}}$, the synchronization condition becomes:
 \begin{equation}\label{sec_bra1}
\frac{1}{\gamma_{1}+\lambda_{1}} \left( \prod_{k=2}^4  \frac{1}{\gamma_{k}+\lambda_{k}} \, +  \prod_{k=5}^7  \frac{1}{\gamma_{k}+\lambda_{k}} \right)  < \sec(\pi/4)^{4}.
 \end{equation}
 If we limit the coupling to the first species only, \refe{sec_bra1} reduces to:
\begin{equation}
\lambda_1  >\frac{\gamma_2\gamma_3 \gamma_4 + \gamma_5\gamma_6 \gamma_7 }{\gamma_2\gamma_3\gamma_4\gamma_5\gamma_6\gamma_7} \cos(\pi/4)^{4} - \gamma_1.
 \end{equation}

 \begin{figure}
 \centerline{
 \includegraphics[scale=0.65]{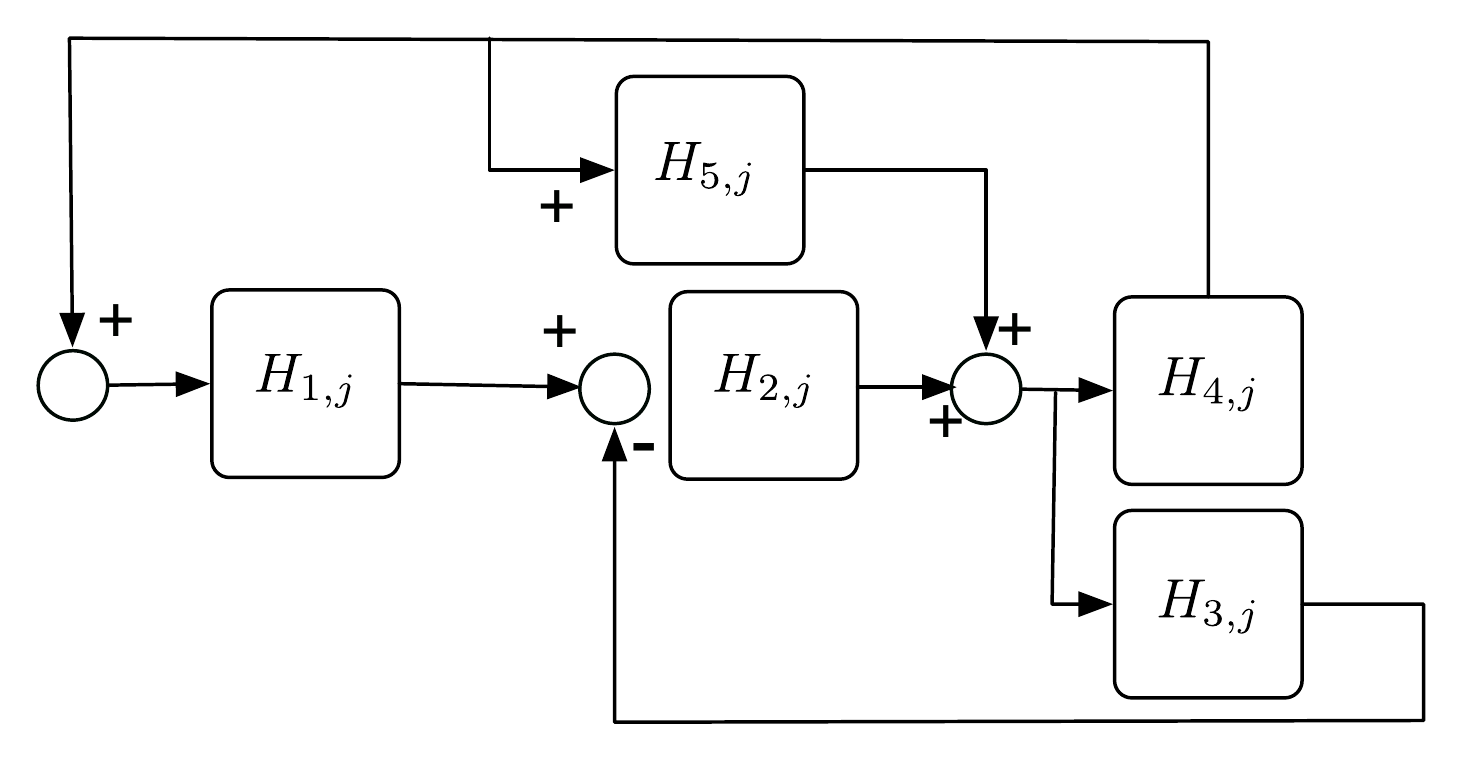}
 }
  \caption{Second type of branched interconnection structure. For notational simplicity the inputs $w_{jk}$ are not shown in the picture.}
  \label{fig:branched2}
 \end{figure}

\noindent  ii) For the interconnection structure depicted in Figure \ref{fig:branched2},  the interconnection matrix is:
\[
 \Sigma_{\footnotesize \mbox{b2}}=
 \bma
 0 & 0 & 0 & 1 & 0\\
 1 & 0 & -1 & 0 & 0\\
 0 & 1 & 0 & 0 & 1\\
0 & 1 & 0 & 0 & 1\\
0 & 0 & 0 & 1 & 0
 \ema,
 \]
and therefore the dissipativity matrix is:
   \[
 E^{\footnotesize \mbox{b2}}_{\tilde \gamma}=
 \bma
 -{\tilde \gamma_{1} } & 0 & 0 & 1 & 0\\
1 & -{\tilde \gamma_{2} } & -1 & 0 & 0\\
 0 & 1 & -{\tilde \gamma_{3} } & 0 & 1\\
 0 & 1 & 0 & -{\tilde \gamma_{4} } & 1\\
 0 & 0 & 0 &1 & -{\tilde \gamma_{5} }\\
 \ema.
 \]
The analysis in \cite{Arcak:2008p782} gives the sufficient condition
\[
\frac{1}{\tilde \gamma_{1} \tilde \gamma_2 \tilde \gamma_{4}} + \frac{1}{\tilde \gamma_4 \tilde \gamma_5} < 1,
\]
which leads to
\begin{equation}\label{branch2}
\frac{1}{\gamma_{4}+\lambda_{4}} \left(  \frac{1}{(\gamma_1+\lambda_{1}) (\gamma_{2}+\lambda_{2})}+  \frac{1}{\gamma_5+\lambda_{5}}\right)  < 1.
\end{equation}
If we limit the coupling to the first species only, \refe{branch2} reduces to:
\begin{equation}
\lambda_1  >\frac{\gamma_5}{\gamma_2(\gamma_4 \gamma_5 - 1)} - \gamma_1.
 \end{equation}
\section{Synchronization in networks of Goodwin oscillators}\label{sec:go}
\begin{figure}
\begin{center}
\begin{tabular}{c c c}
&(a)&\\
&\includegraphics[width = 9.1 cm]{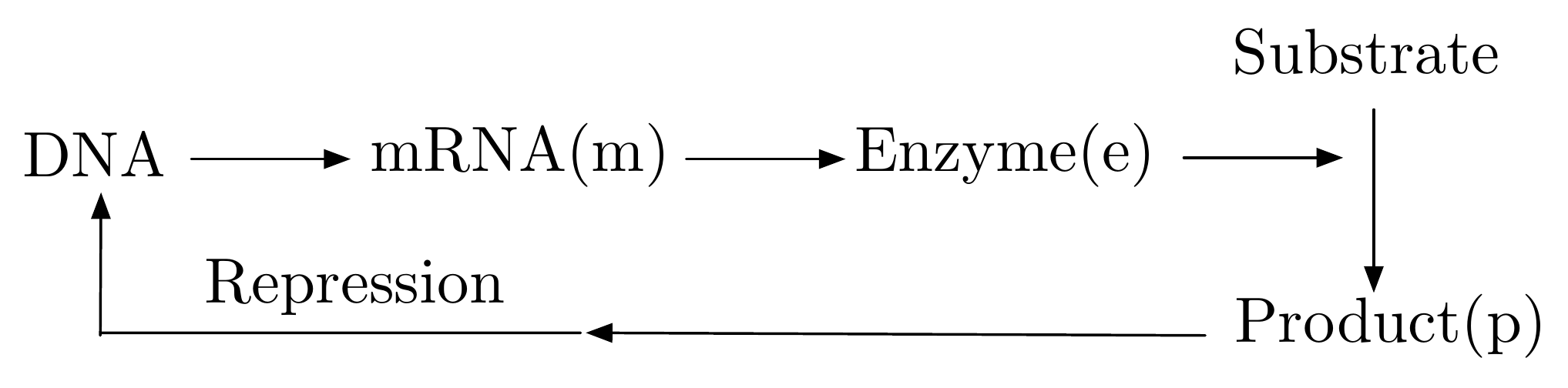}&\\
&&\\
&(b)&\\
&\includegraphics[width = 9.1 cm]{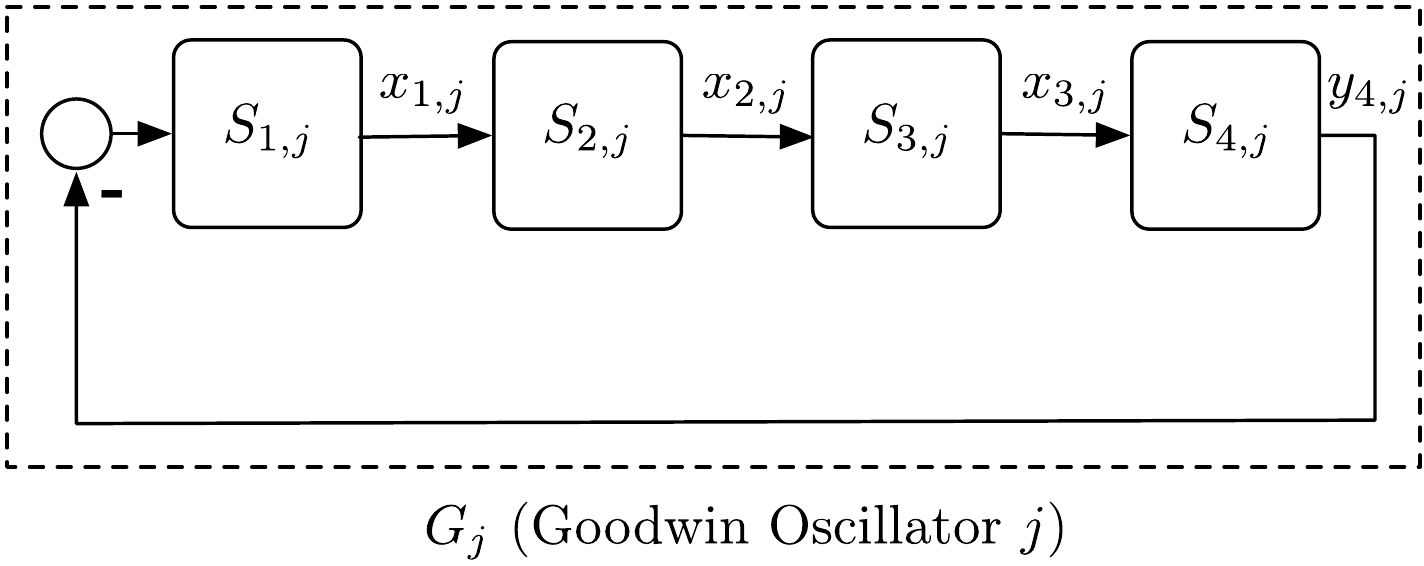}&
\end{tabular}
\end{center}
\caption{(a) Biological interpretation for the Goodwin oscillator. The enzyme (e) combines with the substrate to produce a product (p) which represses the transcription of DNA to mRNA (m), the template for making the enzyme. (b) Input-output scheme representing the mathematical model \refe{goodwin} where $u_{1,j} =  -y_{4,j},\, u_{2,j} =  y_{1,j},\, u_{3,j} =  y_{2,j}$.}
\label{Goodwinscheme}
\end{figure}

We illustrate the proposed theory via a genetic regulatory network example: the Goodwin oscillator. We consider a network of $n$ identical Goodwin oscillators interconnected through a compartmental coupling described by Laplacian matrices $L_{k},\, \kN$. The Goodwin model is an example of cyclic feedback systems described in Section \ref{sec:sec} where metabolites repress the enzymes which are essential for their own synthesis by inhibiting the transcription of the molecule DNA to messenger RNA (mRNA). The model for such a mechanism is schematically shown in Figure \ref{Goodwinscheme}.a and can be described as the cyclic interconnection of $3$ elementary subsystems plus a static nonlinearity (see Figure \ref{Goodwinscheme}.b).

Each Goodwin oscillator can be modeled as a compartment made up of the following four cyclically interconnected sub-systems  (see \cite{FaMaWaTy} for more details):
\begin{equation}\label{goodwin}
\begin{array}{rcll}
S_{k,j}&:&
\left \{
\begin{array}{rcl}
\dot x_{k,j} &=& -b_k \, x_{k,j} + c_k (v_{k,j} + w_{k,j})\\
y_{k,j} &=& x_{k,j}
\end{array} \right., &  \jn, \quad k=1,2,3,  \\
\\
S_{4,j}&:&
\left \{
\begin{array}{rcl}
y_{4,j} &=&\displaystyle -\frac{1}{x_{3,j}^p + 1}
\end{array}
\right., & \jn,
\end{array}
\end{equation}
where $b_k$, $c_k$ are positive coefficients, $p>1$ is the Hill coefficient (that measures the cooperativity of the end product repression) and $w_{k,j}$ are external inputs. The interconnections are encompassed by the inputs
\begin{equation}\label{congod}
v_{k,j} = u_{k,j} + \sum_{z = 1}^n a^k_{j,z} \left(x_{k,z} - x_{k,j}\right), \quad k=1,\ldots,N, \quad j=1,\ldots,n,
\end{equation}
where $u_{1,j} =  -y_{3,j},\, u_{2,j} =  y_{1,j},\, u_{3,j} =  y_{2,j}$, and give rise to the cyclic interconnection matrix
 \[
 \Sigma_{\footnotesize \mbox{g}}=
 \bma
 0 & 0 & 0 & -1\\
 1 & 0 & 0& 0\\
 0 & 1 &0& 0\\
0&0&1&0\\
 \ema,
 \]
 while the second term in \refe{congod} represents the diffusion among the compartments.

From Section~\ref{sec:DiEx}, we observe that the linear sub-systems $S_{k,j}, k=1,2,3$  can be associated to cocoercive operators with constants $\gamma_k = b_k/c_k$. The static nonlinearities $S_{4,j}$ are monotone increasing functions that satisfy $\displaystyle \frac{d}{d\sigma} \left(-\frac{1}{\sigma^p + 1}\right) \leq 1/\gamma_4$, where
\begin{equation}\label{g4}
\gamma_4= \displaystyle \frac{ \left(\sqrt[p-1]{\left(\frac{p-1}{p+1}\right)^p} + 1\right)^2 \left( p+1\right)}{p\left( p-1\right)}, \quad p>1.
\end{equation}
From Section~\ref{sec:DiEx} we know that the co-coercivity coefficient for the static nonlinearity is $\gamma_4$. Since all the blocks are associated to cocoercive operators, Assumption 1 in Theorem 1 is satisfied. The closed loop system is zero-state reachable since it can be fully actuated from the external inputs $w_{k,j}$. Furthermore it is proved in \cite{GuyBartStan:2007p3} that the positive orthant is an invariant set and that the solutions of the closed loop system are bounded. The secant condition for cyclic systems \refe{sec_tot} specifies to
\begin{equation} \label{goodsecgen}
\displaystyle  \left( \gamma_1+\lambda_{1}\right)\left(\gamma_2+\lambda_{2}\right)\left(\gamma_3+\lambda_{3}\right) > c, \quad \quad c=\frac{1}{\gamma_{4} \sec(\pi/4)^{4}}.
\end{equation}
Therefore, if \refe{goodsecgen} holds, then all the conditions of Corollary 1 are satisfied and we conclude that the concentrations of the species in different compartments synchronize when $w_{k,j}=0, k=1,2,3,\; \jn$.

When the compartments are isolated each of them has a unique equilibrium. By choosing the parameters $b_k=c_k=1, k=2,3$ and $b_1 =0.5, c_1 = 1$ it can
be easily proved that the equilibria are asymptotically stable when
$p< 16$. When $p=16$ the steady state undergoes Hopf bifurcation and for
$p=17$ a stable limit cycle arises.
With this choice,  the cocoercive gains are $\gamma_2=\gamma_3 = 1, \gamma_1 = 0.5$. By substituting $p=17$ in \refe{g4} we obtain $\gamma_4 \cong 0.23$ and the secant condition \refe{goodsecgen} becomes
\begin{equation} \label{goodsec}
\displaystyle  \left(0.5+\lambda_{1}\right)\left(1+\lambda_{2}\right)\left(1+\lambda_{3}\right) > c, \quad \quad c=\frac{1}{\gamma_{4} \sec(\pi/4)^{4}} \cong 1.06.
\end{equation}
Condition \refe{goodsec} relates the compartmental coupling to the
synchronization property of the compartments through the algebraic
connectivities $\lambda_k,\, k=1,2,3$.
For simplicity, we assume that those edges that exist all have
the same weight $q$ (which can be interpreted as diffusion coefficients),
i.e. $a^k_{i,j} \in \left\{0,q\right\}$ for every $k=1,2,3$, $i,j =
1,2,\ldots,n$.

Consider for example the case in which only the first and the second species diffuse. Then condition
\refe{goodsec} reduces to
\begin{equation} \label{goodsec2}
\displaystyle  \left(0.5+\lambda_{1}\right)\left(1+\lambda_{2}\right) > c.
\end{equation}
By substituting the expression for the algebraic connectivity (for different graph topologies) in the second order inequality \refe{goodsec2} we can find conditions on the number of cells and the diffusion coefficients such that synchronization is guaranteed. In Table \ref{graphcond} we list these conditions for a number of relevant graph topologies.
\begin{center}
\begin{table}\label{graphcond}
\centerline{
{\footnotesize
\begin{tabular}{|c|c|c|c|}
\hline
Graph &$\lambda_1$&$\lambda_2$& Synchronization condition\\
\hline
$G_c$ (Complete)&&&\\
\multirow{2}{*}{
\includegraphics[width = 2 cm]{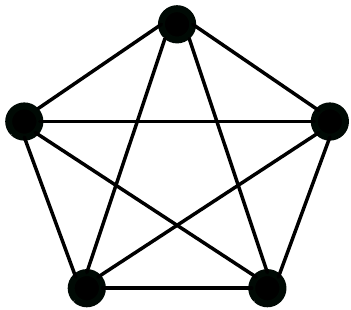} }
& $nq$  & $0$  & $\displaystyle n>\frac{c-0.5}{q}$\\
&&&\\
 \cline{2-4}
 &&&\\
& $nq$  & $nq$ & $\displaystyle n  > \frac{-3+\sqrt{9+8c}}{4q}$ \\
&&&\\
\hline
 $G_s$ (Star)&&&\\
\multirow{2}{*}{
 \includegraphics[width = 2 cm]{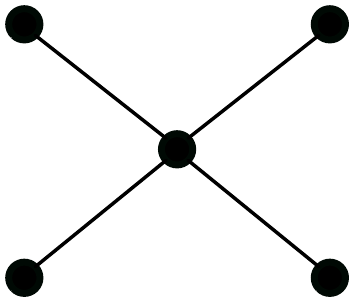}} & $q$ &$0$ & $q>c-0.5$\\
 &&&\\
 \cline{2-4}
 &&&\\
 & $q$ & $q$ & $\displaystyle q  > \frac{-3+\sqrt{9+8c}}{4}$ \\
 &&&\\
\hline
$G_r$ (Ring)&&&\\
\multirow{2}{*}{
\includegraphics[width = 2 cm]{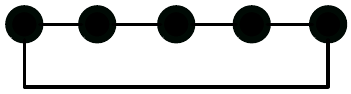} }&  $4 q\sin^2(\frac{\pi}{n})$&$0$&
$
\begin{array}{lll}
q&>&\displaystyle \frac{c-0.5}{4}\\
&&\\
n&<&  \displaystyle  \frac{\pi}{\arcsin\left(\sqrt{ \frac{c-0.5}{4q}}\right)}
\end{array}
$\\
 &&&\\
\cline{2-4}
 &&&\\
&$4 q\sin^2(\frac{\pi}{n})$ &$4 q \sin^2(\frac{\pi}{n})$&
$
\begin{array}{l l l}
\displaystyle q  &>&\displaystyle  \frac{-3+\sqrt{9+8c}}{16}\\
&&\\
{n} &<& \displaystyle  \frac{\pi}{\arcsin\left(\sqrt{ \frac{-3+\sqrt{9+8c}}{16q}}\right)}
\end{array}$\\
 &&&\\
\hline
$G_l$ (Line) &&&\\
\multirow{2}{*}{
\includegraphics[width = 2 cm]{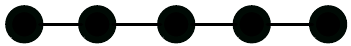} } &$2 q [1 - \cos(\frac{\pi}{n})]$&$0$&
$
\begin{array}{lll}
q&>&\displaystyle \frac{c-0.5}{2}\\
&&\\
n&<&  \displaystyle  \frac{\pi}{\arccos\left( \frac{c-0.5}{2q} +1\right)}
\end{array}
$
\\
 &&&\\
\cline{2-4}
 &&&\\
&$2 q [1 - \cos(\frac{\pi}{n})]$ &$2 q [1 - \cos(\frac{\pi}{n})]$&
$
\begin{array}{l l l}
\displaystyle q  &>&\displaystyle  \frac{\sqrt{9+8c}-3}{8}\\
&&\\
{n} &<& \displaystyle  \frac{\pi}{\arccos\left( \frac{3-\sqrt{9+8c}}{8q}+1\right)}
\end{array}$\\
 &&&\\
\hline
\end{tabular}
}
}
\caption{Sufficient conditions, obtained from \refe{goodsec2}, to achieve synchronization for different diffusive graph topologies.  For each graph depicted in the first column, we consider two cases: (1) only the first species in each compartment are allowed to ``diffuse'', and (2) both the first and second species ``diffuse''.}
\end{table}
\end{center}
\begin{figure}[tbh]\label{complete}
\begin{tabular}{c c}
\includegraphics[width = 7 cm]{{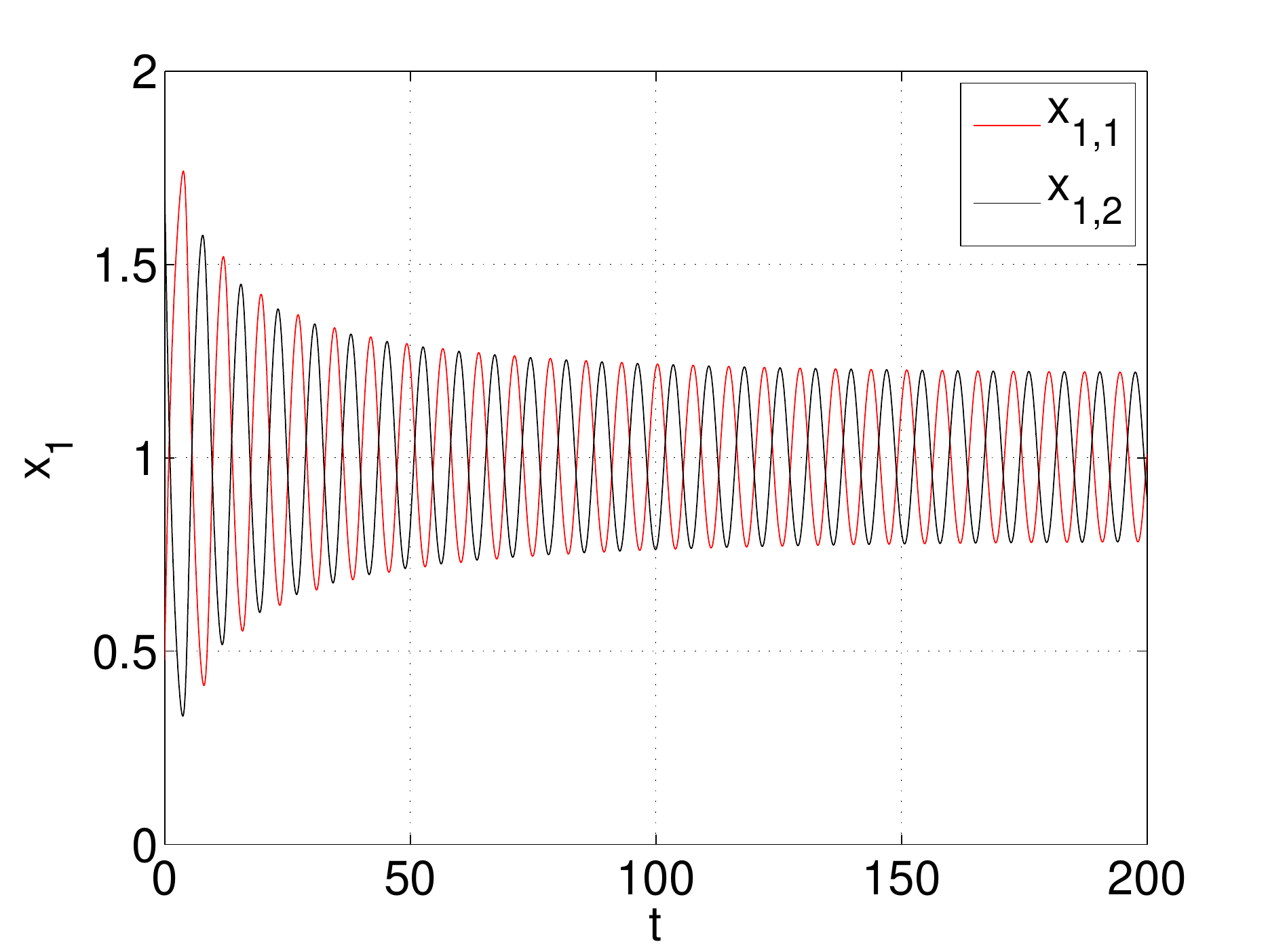}}&\includegraphics[width = 7 cm]{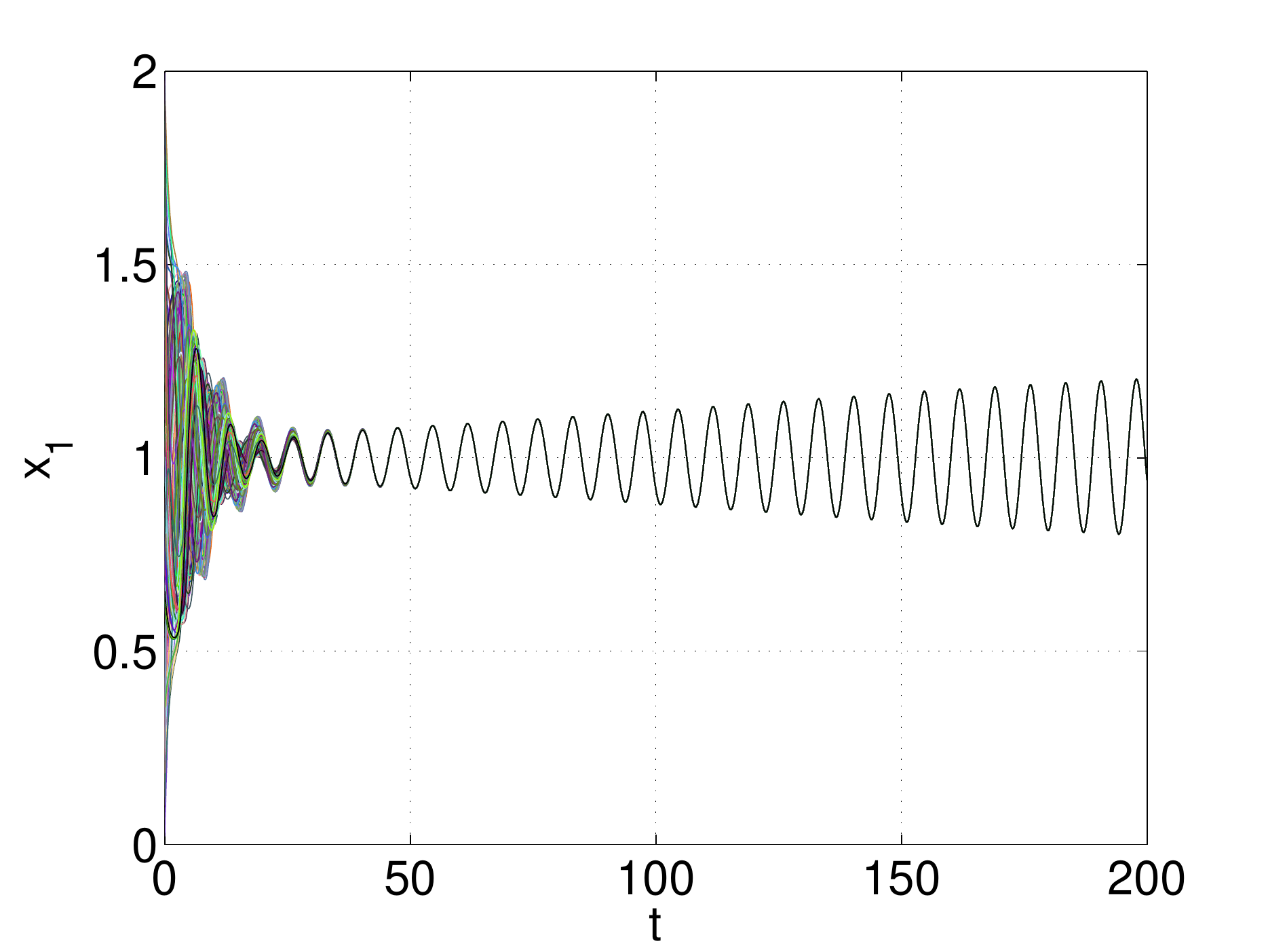}\\
 \includegraphics[width = 7 cm]{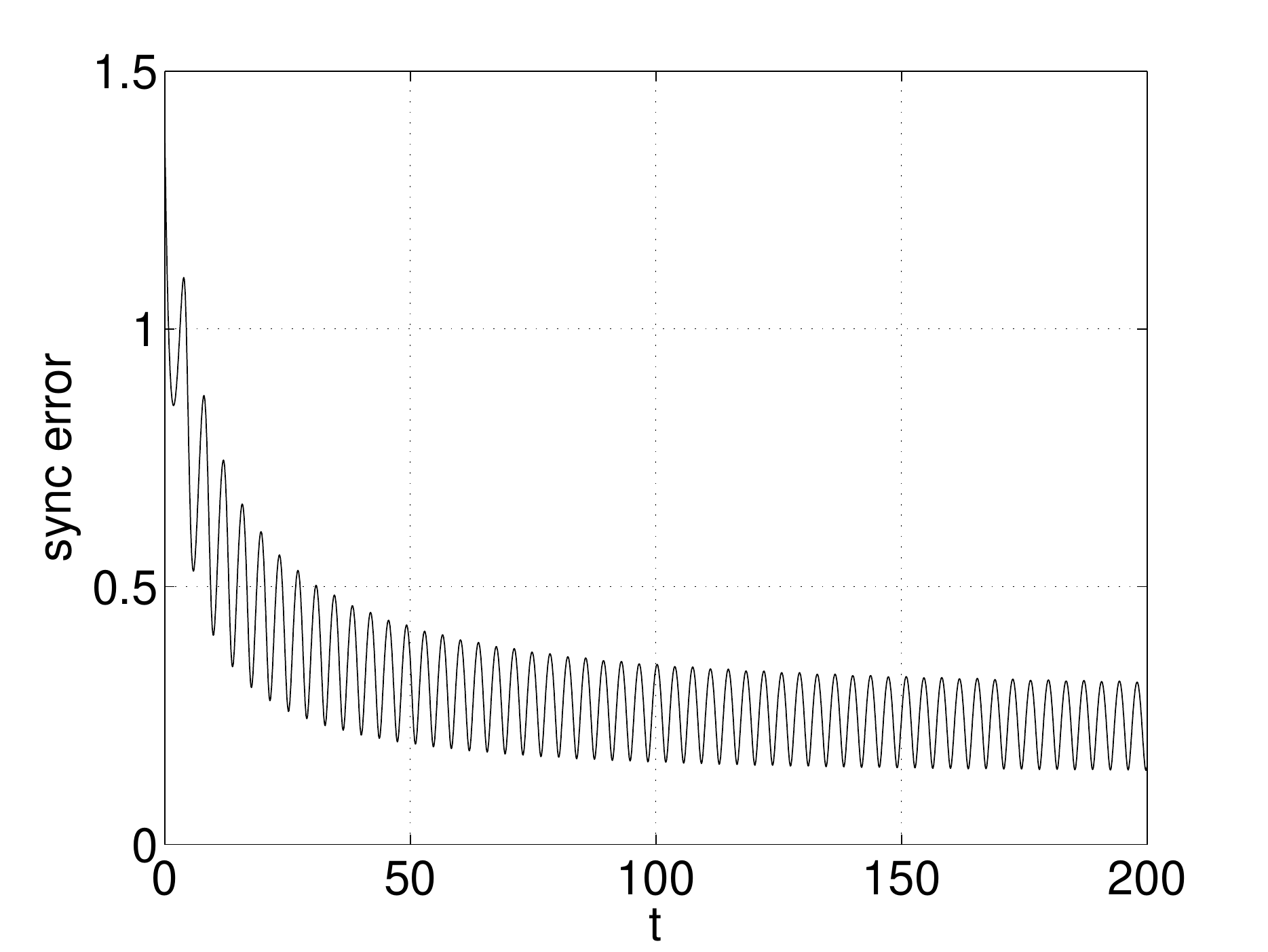}&\includegraphics[width = 7 cm]{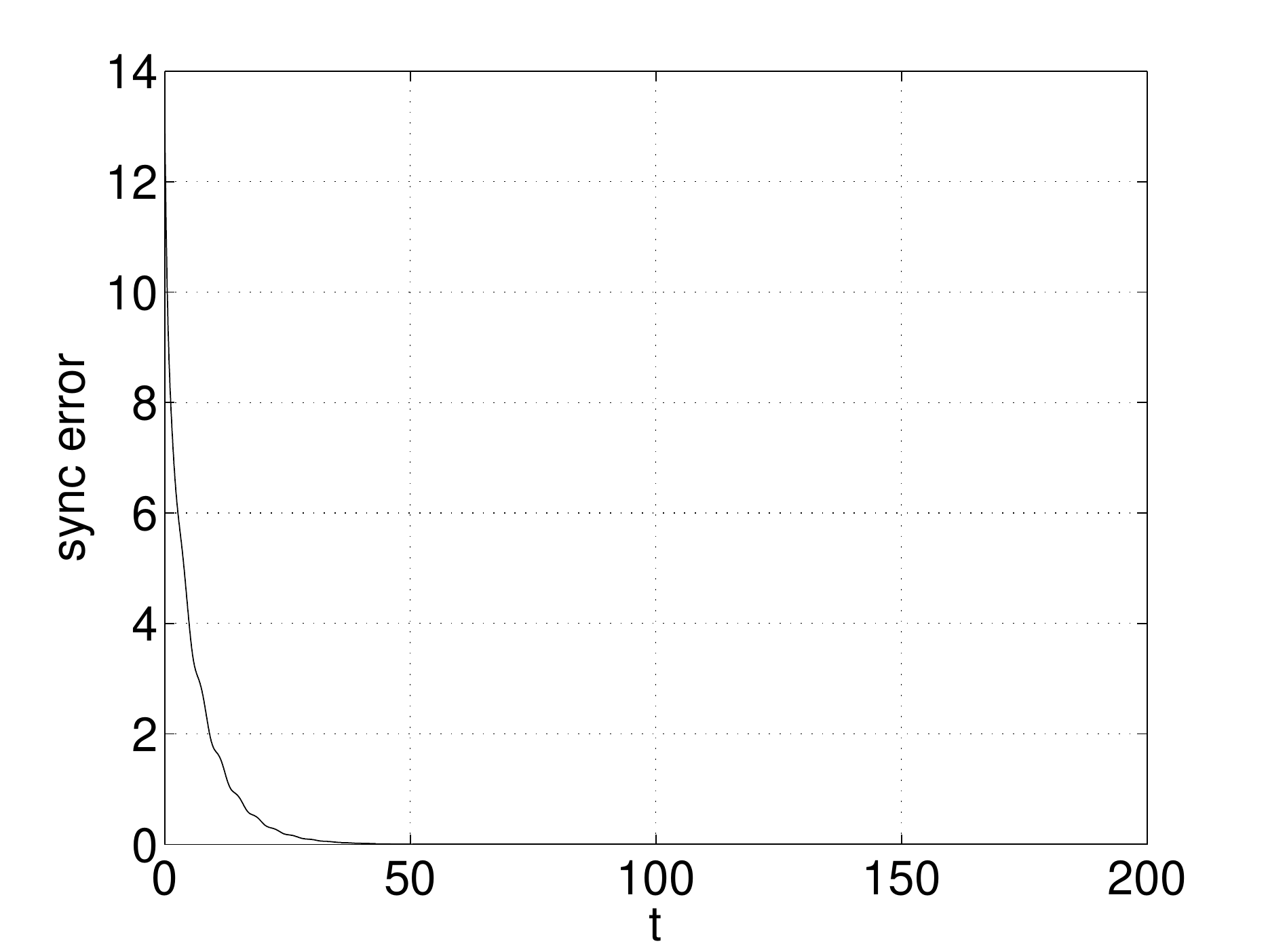}\\
\end{tabular}
\caption{Simulation results for a network of Goodwin oscillators where only the first species are coupled through a complete compartmental coupling $G_c$ and where $q=3\cdot 10^{-3}$. On the left: two oscillators are not sufficient for synchronization. On the right: simulation results for $180$ oscillators. As predicted by the synchronization condition the oscillators synchronize}
\end{figure}
\begin{figure}[tbh]\label{ring}
\begin{tabular}{c c}
\includegraphics[width = 7.1 cm]{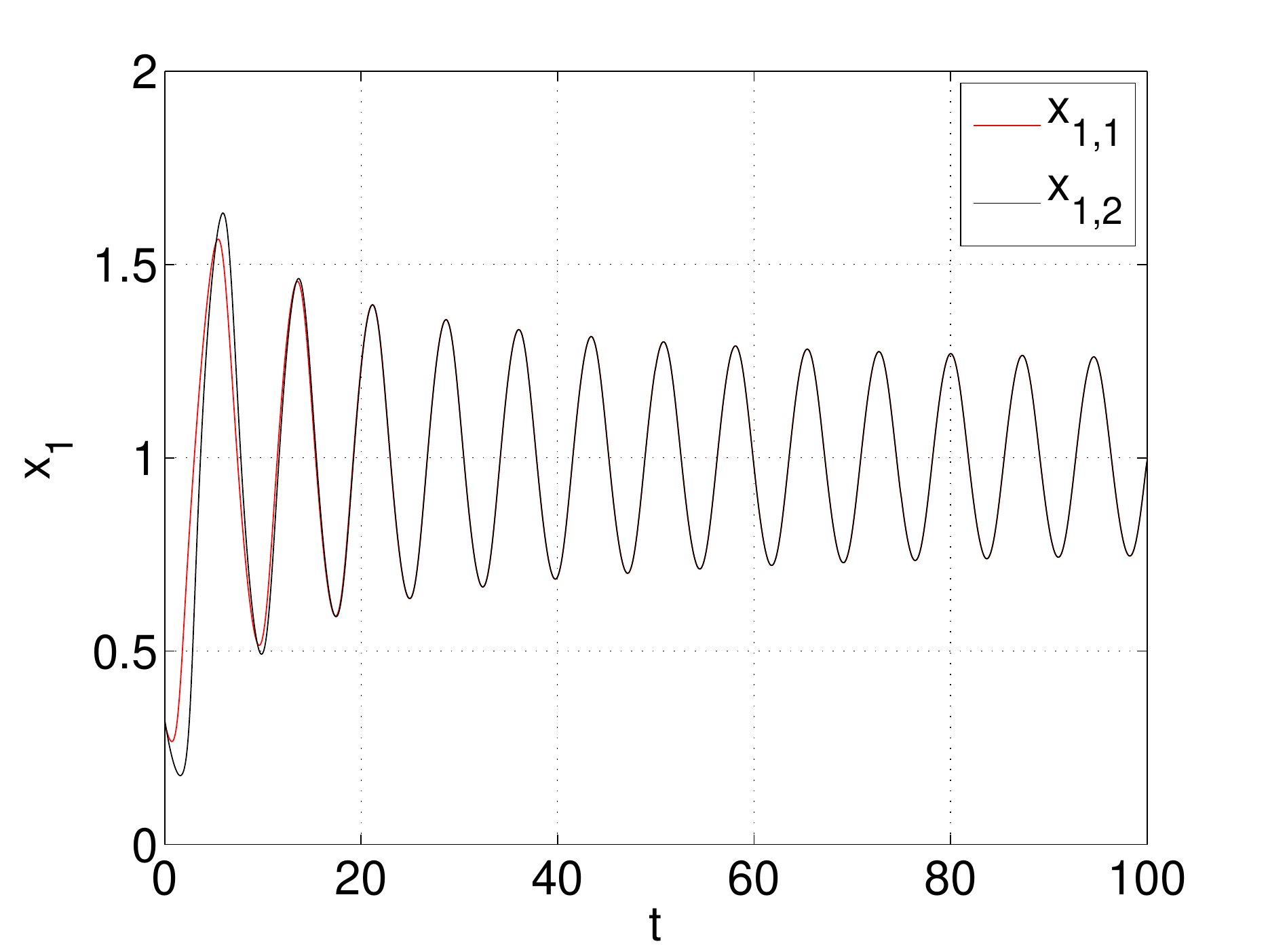}&\includegraphics[width = 7.1 cm]{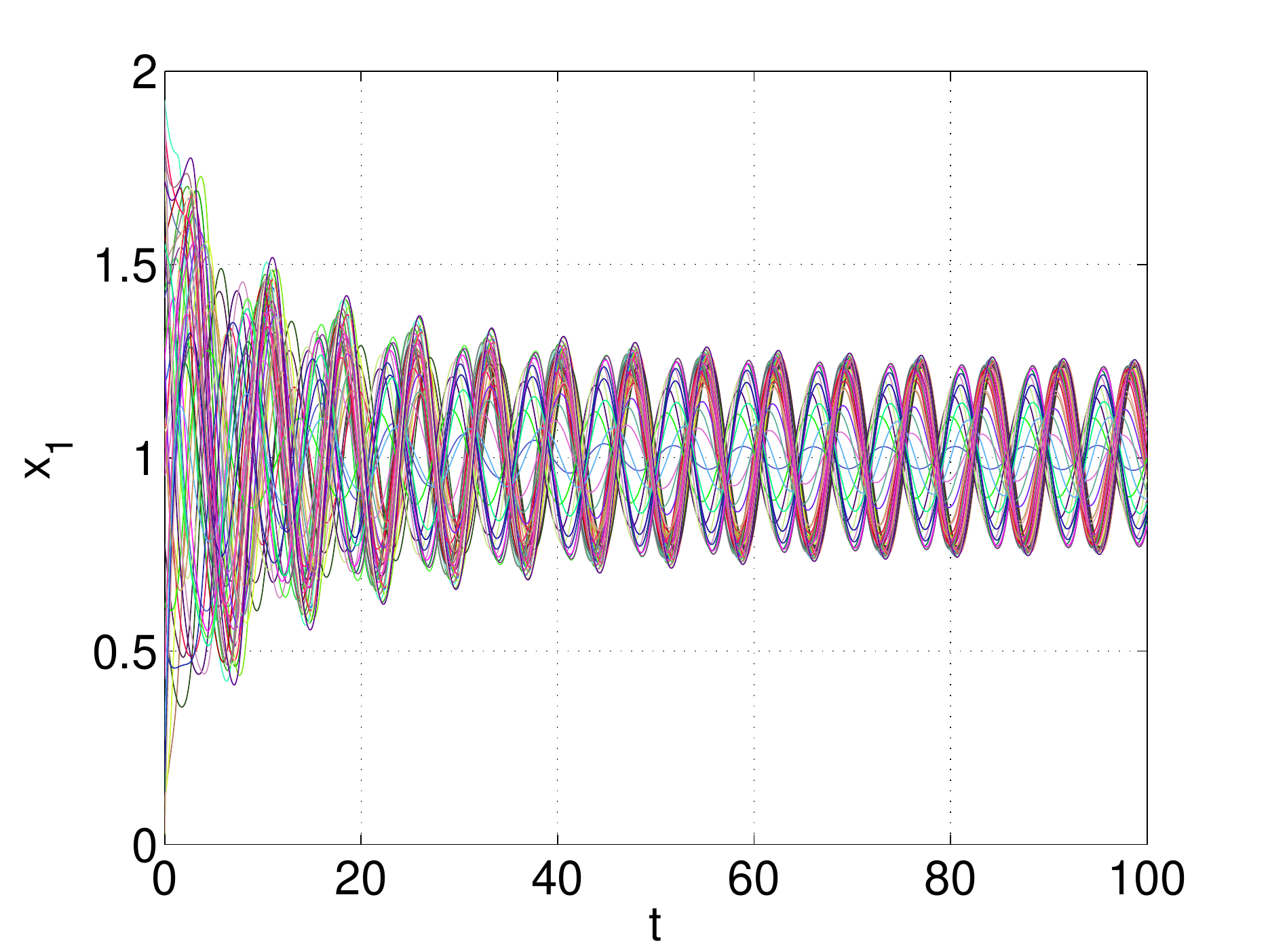}\\
\includegraphics[width = 7.1 cm]{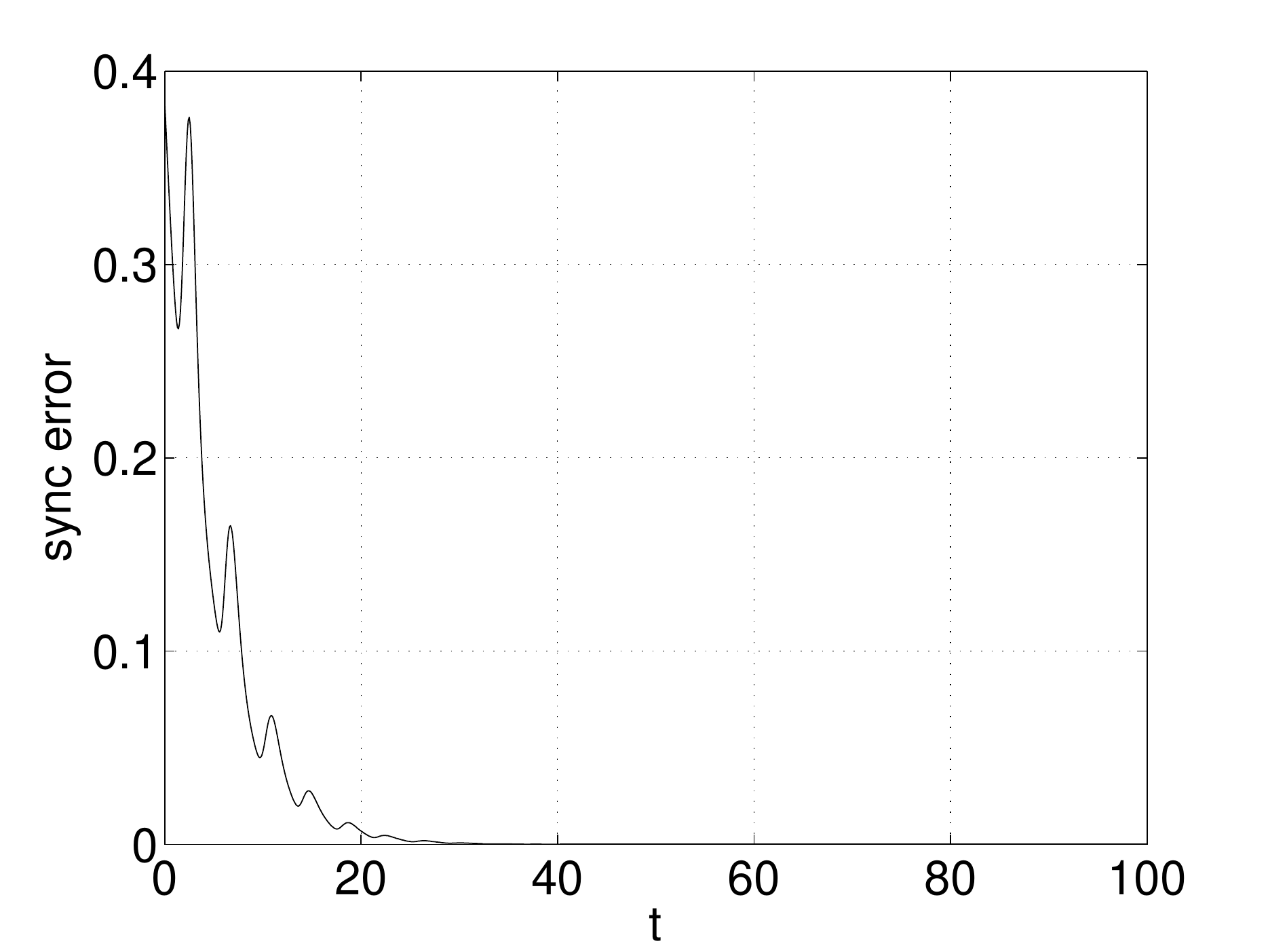} & \includegraphics[width = 7.1 cm]{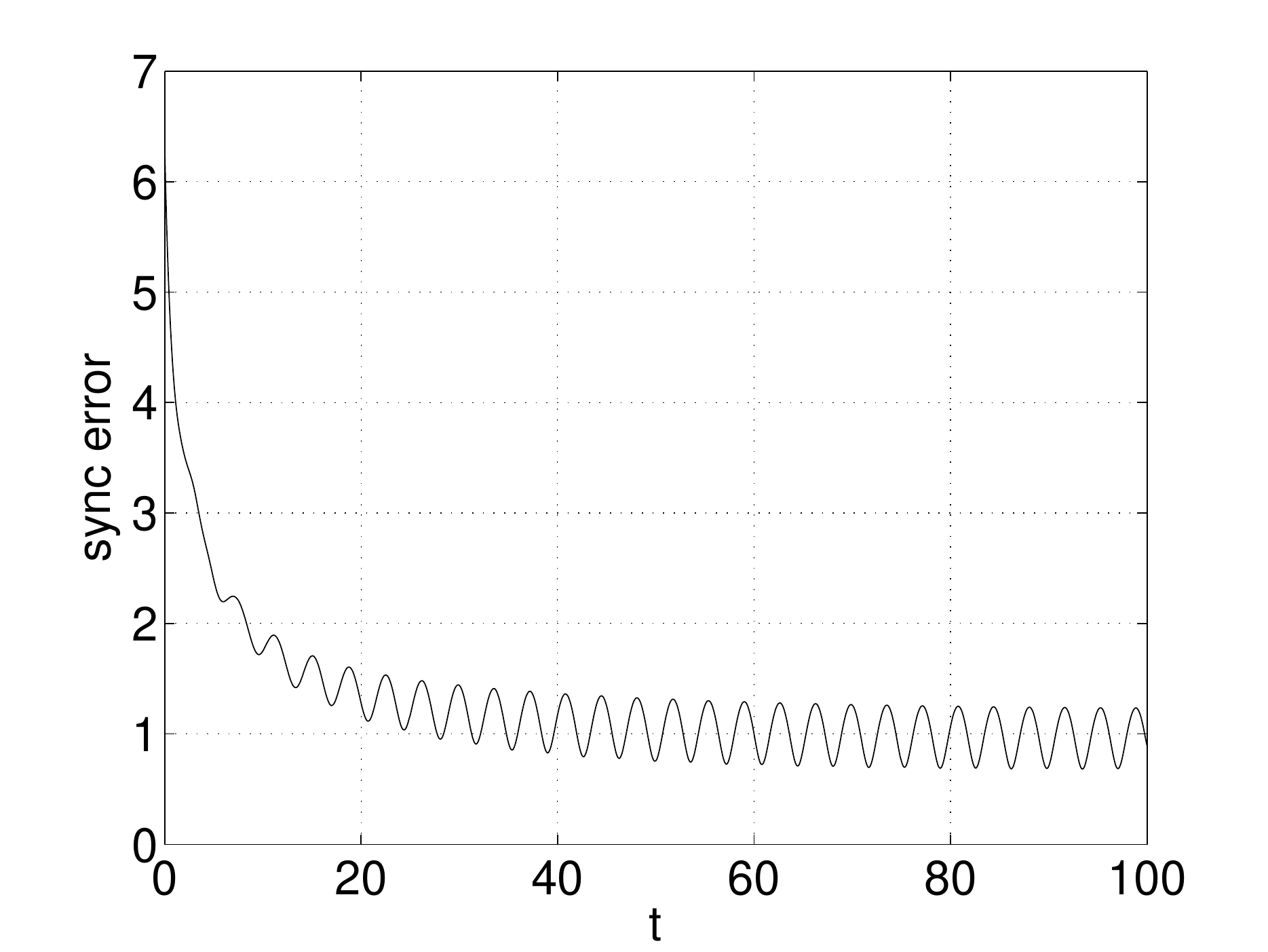}
\end{tabular}
\caption{Simulation results for a network of Goodwin oscillators where only the first species are coupled through a ring compartmental coupling $G_r$ and where $q=0.15$.  On the left: as predicted by the synchronization condition four oscillators are sufficient for the network to synchronize. On the right:  the number of cells is increased up to $45$ and synchronization is not observed.}
\end{figure}

The resulting relations admit interesting biological
interpretations.
Let us think of each compartment as a biochemical network
inside each cell in a population or ``colony'' of $n$ identical cells.

Consider the complete graph  denoted as $G_c$ in Table~\ref{graphcond}.
Now pick a diffusivity coefficient $q$ for which our
``synchronization condition'' estimates fail to hold, and
suppose that the overall network does not synchronize.
From Table~\ref{graphcond}, we observe that a sufficient increase in the
total number of cells in the colony will result in
synchronization between all the cells.
Thus, we may view the number of cells as an order parameter
(or ``synchronization bifurcation'' parameter).
Numerical simulations substantiate this claim, as shown in
Figure~\ref{complete}.

As another example, consider the star graph $G_s$.
Analyzing the conditions in Table~\ref{graphcond}, we see that, for this
type of graph, the number of cells does not play a role in the conditions for
synchronization.  Instead, it is only required now that $q$ be beyond a given
threshold.

The ring graph $G_r$ and the line graph $G_l$ lead to a quite different
qualitative picture.  First, note that we obtain two separate
conditions for $q$ and $n$.
If $q$ is sufficiently large (e.g., in the case of first and second species
coupled, $q>0.074$ for $G_r$ and $q>0.148$ for $G_l$), then we have an upper
bound on the number of cells, instead of a lower bound.
Thus, for either line or ring topologies, in which the graph diameter increases
with the number of cells, we see that a large number of cells, $n$, leads to more restrictive conditions for synchronization.
In Figure~\ref{ring}, this phenomenon is illustrated through simulations.
Indeed, numerous bounds have been derived in the literature \cite{AloMil85,ChuFabMan94} which show that, as the diameter is increased to infinity, the algebraic connectivity of a graph tends to zero.

\subsection{Synchronization conditions and observer design}
To illustrate the idea introduced in Remark 2, we consider two Goodwin oscillators and a directed link coupling the first species. 
\begin{figure}[tbh]
\begin{center}
\includegraphics[width = 10.1 cm]{{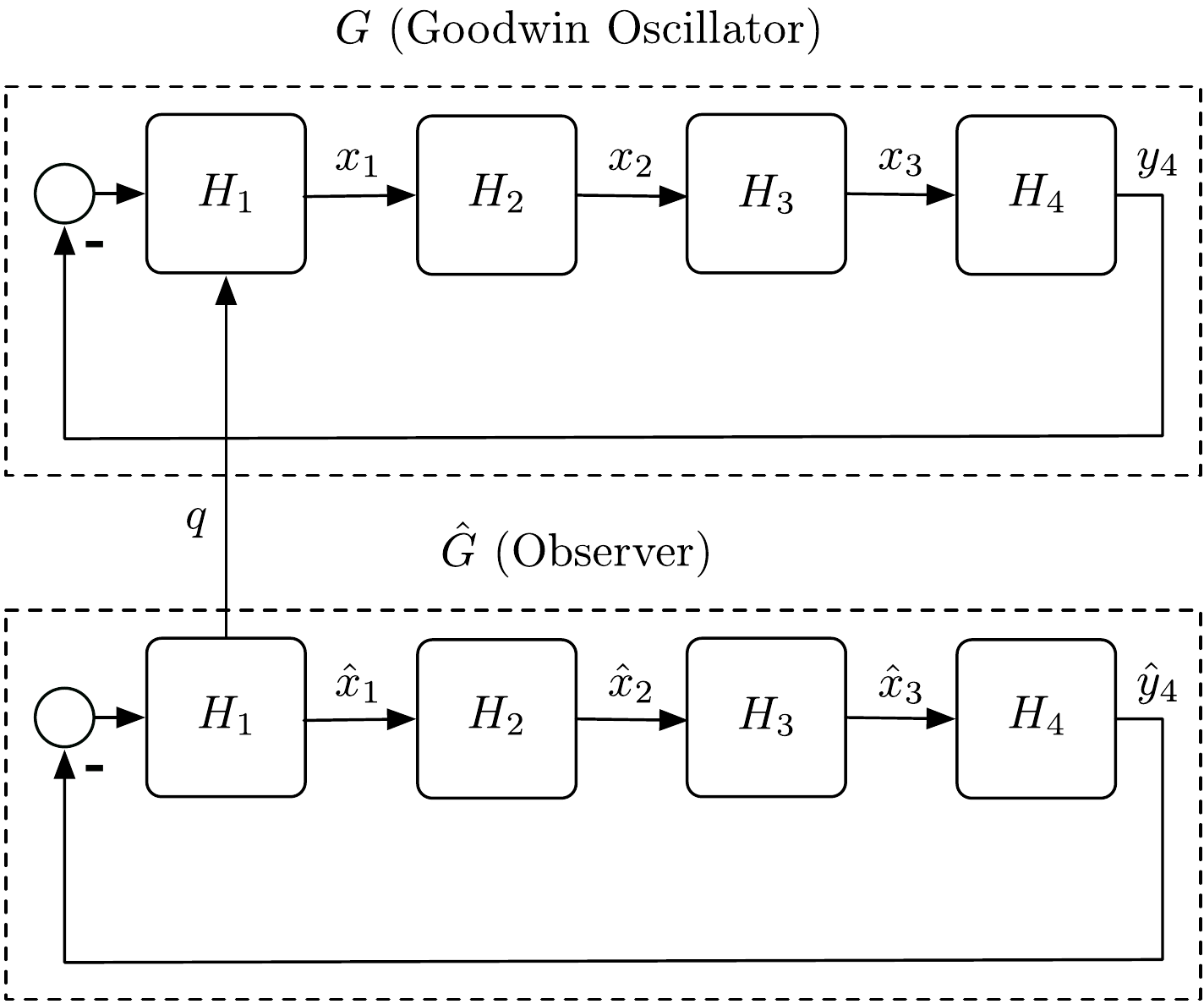}}
\end{center}
\caption{System-observer interpretation for two Goodwin oscillators unidirectionally coupled}
\label{fig:good_obs}
\end{figure}
This special interconnection structure gives rise the following system-observer dynamics:
\begin{equation*}\label{good_obs}
\begin{array}{rcl}
G&:&
\left \{
\begin{array}{rcl}
\dot x_{1} &=& -0.5 \, x_{1} - y_{4} \\
\dot x_{2} &=& -0.5 \, x_{2} + 0.5 x_{1}\\
\dot x_{3} &=& -0.5 \, x_{3} + 0.5 x_{2}\\
y_{4} &=& \displaystyle -\frac{1}{x_{3}^p + 1}
\end{array} 
\right.\\
\\
\hat G&:&
\left \{
\begin{array}{rcl}
\dot {\hat x}_{1} &=& -0.5 \, \hat x_{1} - \hat y_{4} + q(x_{1} - \hat x_1) \\
\dot {\hat x}_{2} &=& -0.5 \, \hat x_{2} + 0.5 \,\hat x_{1}\\
\dot {\hat x}_{3} &=& -0.5 \, \hat x_{3} + 0.5 \,\hat x_{2}\\
\hat y_{4} &=& \displaystyle -\frac{1}{\hat  x_{3}^p + 1}
\end{array} \right.
\end{array}
\end{equation*}
The term $q(x_{1} - \hat x_1)$ (where $q$ is the weight of the link), that was interpreted as diffusion of the first species concentrations, is now the output injection to the observer $\hat G$ (see Figure~\ref{fig:good_obs}).
Then the synchronization condition \refe{goodsec} reduces to
\begin{equation} \label{goodsec_obs}
\displaystyle  0.5+2\,q  > c, \quad \quad c \cong 1.06,
\end{equation}
and can be interpreted as a sufficient condition for the observer error to converge to zero.
We conclude that if
\[
q>\displaystyle \frac{c - 0.5}{2},
\]
then the errors $x_k - \hat x_k \rightarrow 0$ as $t\rightarrow \infty$, $k=1,2,3$.

\section{Conclusion and future work}
Synchronization properties for networks of nonlinear
systems have been investigated combining the input-output properties of the 
subsystems with the information about the structure of network. The proposed model is motivated by cellular networks where signaling occurs both internally, through interactions of species, and externally, through intercellular signaling. 
Results for state-space models as well as biochemical applications have been derived as corollaries of the main result. The extension of the present work to diffusion models (by using partial differential operators) is currently being developed by the authors. 

\bibliography{arXiv_08Mar09.bib}
\end{document}